\documentclass[smallextended]{amsart}
\usepackage[table,usenames,dvipsnames]{xcolor}
\usepackage{amsmath}
\usepackage{amssymb}
\usepackage{graphicx}

\usepackage{epstopdf}
\usepackage{hyperref}
\usepackage[dvipsnames]{xcolor}
\usepackage{enumitem}
\usepackage{bm}
\usepackage{tikz}
\usepackage{siunitx}
\usetikzlibrary{positioning,arrows}
\setlist{nolistsep}

\newcommand{\epsi}{\ensuremath{\varepsilon}}
\newcommand{\Sym}{\mathrm{Sym}}
\newcommand{\aSym}{\mathrm{aSym}}

\newcommand{\grad}{\nabla}%
\newcommand{\dive}{\grad\cdot}%
\newcommand{\ddiv}{\texttt{div}}%
\newcommand{\dgrad}{\texttt{grad}}%

\DeclareMathOperator{\trace}{tr}
\DeclareMathOperator{\inv}{inv}
\newcommand{\Fint}{F_{\text{int}}}
\newcommand{\Real}{\mathbb{R}}
\newcommand{\Pol}{\mathbb{P}}
\newcommand{\proj}{\mat{\pi}}
\newcommand{\dproj}{\mat{P}}
\newcommand{\mat}[1]{{\boldsymbol{#1}}}
\newcommand{\vect}[1]{{\boldsymbol{#1}}}
\newcommand{\dx}{\,\text{d}\vect{x}}
\newcommand{\brac}[1]{\left<#1\right>}
\newcommand{\id}{\mat{I}}
\newcommand{\Pcal}{\mathcal{P}}
\newcommand{\bigand}{\quad\text{ and }\quad}
\newcommand{\Vcal}{\mathcal{V}}
\newcommand{\Tcal}{\mathcal{T}}
\newcommand{\norm}[1]{\left\|#1\right\|}
\newcommand{\abs}[1]{\left|#1\right|}
\newcommand{\projgrad}{\proj^{\grad}\hspace{-0.6mm}}
\newcommand{\projgradi}{\proj_i^{\grad}\hspace{-0.6mm}}
\newcommand{\quot}{/}
\newcommand{\smallinvroottwo}{{\small\frac{1}{\sqrt{2}}}}
\newcommand{\smallmininvroottwo}{{\small\frac{-1}{\sqrt{2}}}}

\def\co2{$\text{CO}_2$}
\newcommand{\trans}{T}
\newcommand{\insertcite}[1]{\textcolor{red}{[insert reference\ifx\empty#1\else\ to #1\fi]} }
\newcommand{\fracpar}[2]{\frac{\partial #1}{\partial #2}}
\newcommand{\sbed}{\textit{sbed} }
\let\vphi=\varphi
\newcommand{\smallheading}[1]{{\it\bf #1}}

\usepackage{soul}
\soulregister\cite7
\soulregister\ref7
\soulregister\eqref7
\soulregister\num7
\soulregister\pageref7
\sethlcolor{yellow}

\title{Virtual Element Method for Geomechanical Simulations of Reservoir Models}

\author{Odd Andersen \and Halvor M. Nilsen \and Xavier Raynaud}

\begin{document}

\maketitle
\begin{abstract}
  In this paper we study the use of Virtual Element Method (VEM) for
  geomechanics. Our emphasis is on applications to reservoir simulations. The
  physical processes behind the formation of the reservoirs, such as sedimentation,
  erosion and faulting, lead to complex geometrical structures. A minimal
  representation, with respect to the physical parameters of the system, then
  naturally leads to general polyhedral grids.  Numerical methods which can directly
  handle this representation will be highly favorable, in particular in the setting
  of advanced work-flows. The virtual element method is a promising candidate to
  solve the linear elasticity equations on such models. In this paper, we investigate
  some of the limits of the VEM method when used on reservoir models. First, we
  demonstrate that care must be taken to make the method robust for highly elongated
  cells, which is common in these applications, and show the importance of
  calculating forces in terms of traction on the boundary of the elements for
  elongated distorted cells. Second, we study the effect of triangulations on the
  surfaces of curved faces, which also naturally occur in subsurface models. We also
  demonstrate how a more stable stabilization term for reservoir application can be
  derived.
\end{abstract}
\newpage

\section{Introduction}

\begin{figure}[h]
  \centering
  \includegraphics[width=0.7\textwidth]{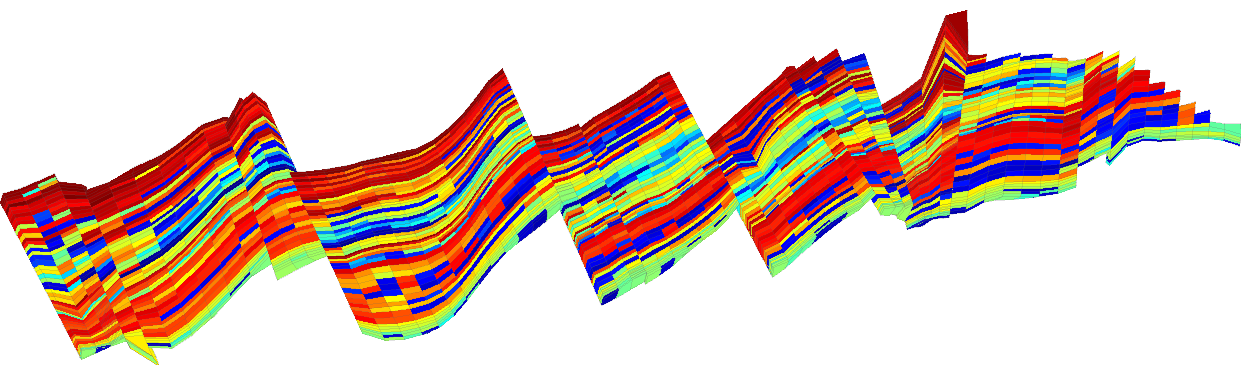}
  \caption{Section of the Gullfaks reservoir model (Norway). Each color represents a
    different material property. We observe how the large aspect ratio in the cells.}
  \label{fig:gullfakssection}
\end{figure}

Sedimentary formations are the result of long and complex geological
processes. Sedimentation creates thin layers, faulting creates nontrivial connections
between the layers and erosion creates degenerate layers. The formation retains an
overall stratigraphic structure, in the sense that very different spatial
correlations in the material properties can be observed between the horizontal and
vertical directions, and long and thin cells are specific to reservoir simulations,
see the section represent in Figure \ref{fig:gullfakssection}.
The geometric modeling of sedimentary formations requires the parameterization of a
very large number of complicated interfaces. Each interface then separates regions
with material properties that may differ of several order of magnitude and must be
captured with maximum accuracy. Because of these difficulties, computational
considerations are often not prioritized in the design of geological grids, which
will typically contain highly irregular cell shapes. The grid and material properties
are strongly related, which cause severe limitations on remeshing. The industry
standard for reservoir grids is the corner-point format \cite{ponting1989corner}. In
a corner-point grid, pillars which have a dominant vertical direction are first
defined from a two-dimensional Cartesian partition. Then, for each set of four
adjacent pillars, hexahedron cells are constructed by choosing 2 points on each
pillars and connecting all these points (see the detailed in Section
\ref{subsec:compaction3D}). Many geometrical grid formats have been proposed to
improve on this format, for example Skua Grid \cite{Gringarten2008}, S-Grid, Faulted
S-Grid and Cut-Cell \cite{Mallison2014}. By refining the mesh, it is of course always
possible to improve the quality of the mesh from the point of view of numerical
computation, but all compact representation of the underlying geology, that is a
representation where the data (the material properties) is represented by the minimum
number of cells, will lead to cells with high aspect ratio, distorted cells, faces or
cells of very different sizes, cells or faces with different shapes.  Methods which
are robust for such grid will greatly simplify the modeling of subsurface physics.

In recent years, the coupling of geomechanical effects with subsurface flow has
become more and more important in many areas including: oil and gas production from
mature fields, oil and gas production from fractured tight reservoirs, fractured rock
for geothermal application and risk assessment of \co2 injection. Realistic modeling
of these applications is hampered by the differences in the way geomechanics and flow
models are build and discretized. Traditionally, the mechanic problems are solved
using finite element methods but they are difficult to adapt to the standard
geometrical representation of reservoir models, such as corner-point grids. In
contrast, the Virtual Element Method (VEM) can operate on general polyhedral
grids. As such, the ability of the method to handle irregular grids makes it very
attractive for geomechanical applications. In this paper, we investigate if the
method can effectively be applied on realistic reservoir grids. Our main result
concerns the treatment of the load term. We observe that standard stabilization terms
presented in the literature are not adapted to elongated cells with large aspect
ratios, which are standard in geological models. We propose a modification of the
stabilization constant which can be used in the 2D case and a discrete gradient
approach to compute the load term which turns out to be little sensitive to the
choice of stabilization and can be easily extended to 3D. In a first part of this
paper, we present the VEM method following mostly \cite{gain2014} but we also try to
clarify the connection with the construction of the projection operators, as
introduced in the basic principles of \cite{beirao2013basic}. In our numerical
experiments, we will focus on the performance of the VEM method on geological grids.
The emphasis will be on corner-point grids and complex small scale sedimentary
models. To be able to demonstrate the VEM method in the reservoir setting we used the
MRST framework \cite{Lie12:comg} to simplify the grid handling. The numerical
implementation of the VEM method used in this paper can be downloaded from
\cite{MRST:2016}, in particular the test case concerning the Norne reservoir model
(see Section \ref{subsec:compaction3D}) is readily available from there.

\section{The equations of linear elasticity}

We consider the equation of linear elasticity for small deformations. The
displacement is given by $\vect{u}(\vect{x})\in\Real^d$ for $d = 2, 3$ (2D or 3D
case) and $\vect{x}\in\Omega\subset\Real^{d}$. The equations are given by
\begin{equation}
  \label{eq:lin_elast_cont}
  \dive\mat\sigma = \vect{f}
\end{equation}
with
\begin{equation}
  \label{eq:defsigmaepsi}
  \mat\sigma =  \mat{C}\mat\epsi\quad\text{ and }\quad \mat\epsi = \frac{1}{2}(\nabla +\nabla^\trans)\vect u
\end{equation}
Here, $\mat\sigma$, $\mat\epsi$ and $\vect{u}$ denote the Cauchy stress tensor, the
infinitesimal strain tensors and the displacement field, respectively. The vector
function $\vect{f}:\Real^d\to\Real^d$ is an external volumetric force that we will
refer to as the \textit{load term}. The linear operator $\mat{C}$ is a fourth order
stiffness tensor which satisfies, for some constant $c>0$, the ellipticity condition
\begin{equation}
  \label{eq:elliptcond}
  c\,\mat{S}:\mat{S} \leq \mat{S}:\mat{C}\mat{S},
\end{equation}
for any symmetric matrix $\mat{S}\in\Real^{d\times d}$. The symbol $:$ denotes the
scalar product in $\Real^{d\times d}$ defined as
\begin{equation}
  \label{eq:scalprodRdd}
  \mat{A} : \mat{B} = \trace(\mat{A}^\trans\mat{B}),
\end{equation}
for any $\mat{A},\mat{B}\in\Real^{d\times d}$. 

\section{Presentation of the VEM method for linear elasticity}

The VEM method was first introduced in the framework of mimetic discretization
methods but later rephrased in the language of finite element methods (see
\cite{da2014mimetic} for discussions). A general presentation of VEM is given in
\cite{beirao2013basic}. The same authors present convergence results for linear
elasticity in \cite{da2013virtual}. Practical details on the implementation of VEM
are given in \cite{beirao2014hitchhiker}. Our implementation of the VEM follows the
presentation done in \cite{gain2014} where the specific case of linear elasticity is
considered. We rewrite the equation of linear elasticity in the weak form,
\begin{equation}
  \int_{\Omega}\mat\epsi(\vect{v}) : \mat{C} \mat\epsi (\vect{u}) \dx = \int_{\Omega} \vect{f}\cdot\vect{v}\dx,
\end{equation}
which must hold for all displacement field $\vect{v}:\Real^d\to\Real^d$. Let $N_c$
denote the number of cells and $\left\{E_i\right\}_{i=1}^{N_c}$ the grid cells. We
define the bilinear form $a_{E_i}$ as
\begin{equation*}
 a_{E_i}(\vect{u},\vect{v}) = \int_{E_i} \mat\epsi(\vect{u}) : \mat{C} \mat\epsi(\vect{v})\,\dx 
\end{equation*}
and decompose the global bilinear energy form $a(\cdot,\cdot)$ in cell contributions,
\begin{equation}
  \label{eq:abil}
  a(\vect{u},\vect{v}) := \int_\Omega \mat\epsi(\vect{u}) : \mat{C} \mat\epsi(\vect{v})\,\dx = \sum_{i=1}^{N_c} a_{E_i}(\vect{u},\vect{v}).
\end{equation}
In the rest of this section, we consider a given cell $E$ and will denote by
$\Vcal_E$ the finite dimensional approximating function space in $E$. In the VEM
approach, the basis functions of $\Vcal_E$ are not known explicitly but, for a
first-order VEM method, the requirements on $\Vcal_E$ are that it contains the space
of polynomials of order 1, denoted $\Pol_1(E)$, and that the bilinear form
$a_E(\vect{u},\vect{v})$ can be computed exactly for any $\vect{u}\in\Pol_1(E)$ and
any $\vect{v}\in\Vcal_E$, using only the degrees of freedom of $\vect{v}$. As in the
standard finite element method, the degrees of freedom are the nodal displacements,
so that the continuity at the boundaries of each element is ensured by requiring
linearity on the edges and a local reconstruction on the faces, which depends only on
the values at the edges of the face where the reconstruction is done. The system
matrix can be assembled element-wise. Let us denote
\begin{equation}
  \Vcal_E^{\text{scalar}} = \{\vect{v}\in H^1(E)\ |\ \vect{v}_{|e} \in \Pol_1(e)\text{ for all edges }e\},
\end{equation}
for $i=1,\ldots,d$. For a given node $\eta$ of $E$, we can construct a function
$\phi_\eta$ in $\Vcal_E^{\text{scalar}}$ such that $\phi_{\eta}(\bar\eta) = 1$ if
$\bar\eta = \eta$ and zero if $\bar\eta\neq\eta$. The virtual basis functions of
$\Vcal_E$ are then given by
\begin{equation}
  \label{eq:defbasis}
  \vect{\phi}_{\eta}^{k}(\vect{x})=\phi_\eta(\vect{x})\vect{e}_k 
\end{equation}
for $\eta\in N(E)$ and $k=1,\ldots,d$, where $N(E)$ denotes the set of nodes of the
cell $E$ and $e_k$ is the unit vector in the direction given by the index $k$. After
having introduced the projection operator, we will add some extra requirements for
$\phi_\eta$ concerning its first and second order moment. But beside that, no more
explicit properties for $\phi_\eta$ are needed and this is one of the important point
of the method, which also makes it so flexible. The projection operator, which we
denote $\projgrad$, is defined with respect to the energy norm $a_E$. We consider
first order approximations and, for any displacement field $\vect{u}$ in the Hilbert
space $[H^1(E)]^3$, the projection $\projgrad(\vect{u})$ of $\vect{u}$ is defined as
the element $\vect{p}\in[\Pol_1(E)]^3$ such that
\begin{equation*}
  a_E(\vect{p},\vect{q}) = a_E(\vect{\vect{u}},\vect{q}),
\end{equation*}
for all $\vect{q}\in[\Pol_1(E)]^3$. Since the bilinear form $a_E$ is degenerate,
additional conditions must be imposed to define completely $\projgrad$, see
\eqref{eq:defprojgrad} in Section \ref{subsec:projoperator} for the rigorous
definition. For any displacement field $\vect{u}$, the energy
$a_E(\vect{u},\vect{u})$ can be decomposed using Pythagoras' identity,
\begin{equation}
  \label{eq:pythagoras}
  a_E(\vect{u},\vect{u}) = a_E(\projgrad\vect{u},\projgrad\vect{u}) + a_E((\id - \projgrad)\vect{u},(\id - \projgrad)\vect{u}).
\end{equation}
The first term on the right-hand side can be computed exactly from the degree of
freedom, for any $\vect{u}\in\Vcal_E$. The last term can not be handled generically
and is therefore replace by a \textit{stabilization} term which takes the form of a
bilinear form $s_E$, whose role is to ensure that the ellipticity of $a_E$ is
retained. Hence, the energy is finally approximated by
\begin{equation}
  \label{eq:pythagoras}
  a_{h,E}(\vect{u},\vect{u}) = a_E(\projgrad\vect{u},\projgrad\vect{u}) + s_E((\id - \projgrad)\vect{u},(\id - \projgrad)\vect{u}).
\end{equation}
In the general framework of VEM, as introduced in \cite{beirao2013basic}, the
computation of the projection operator typically requires the computation of an
inverse, locally for each cell. The formulation in \cite{gain2014} has the advantage
of giving an explicit expression of the projection operator. In the presentation that
follows, we will try to clarify the connection between the two approaches.

\subsection{The kinematics of affine displacement}

The physics of linear elasticity is associated with linear deformations, in
particular the rigid body motions play a crucial role. Let us recall some simple
facts on the kinematics of affine displacements. The linear space of affine
displacements, which we denote by $\Pcal$, corresponds to the sum of the translations
and linear transformation so that any $\vect{l}\in \Pcal$ can be written as
$\vect{l}(\vect{x}) = \vect{u} + \mat{L}\vect{x}$, for $\vect{u}\in\Real^d$ and
$\mat{L}\in\Real^{d\times d}$. The dimension of the space of $\Pcal$ is $d^2 +
d$. The subspace of rigid body motion, which we denote $\Pcal_r$, contains the
rotation and the translation. Any $\vect{l}\in\Pcal_r$ can be written as
\begin{equation}
  \label{eq:expPr}
  \vect{l}(\vect{x}) = \vect{u} + \mat{\Omega}(\vect{x} - \vect{x}_0),
\end{equation}
for any $\vect{u},\ \vect{x}_0\in\Real^{d}$ and $\mat{\Omega}$ that belongs to the
space of skew-symmetric matrices, denoted $\aSym(\Real^d)$. There is a redundancy in
the choice of $\vect{x}_0$ and $\vect{u}$ so that a unique decomposition of
$\vect{l}\in\Pcal_r$ is given by $\vect{l}(\vect{x}) = \vect{u} +
\mat{\Omega}\vect{x}$ for $\vect{u}\in\Real^d$ and $\mat{\Omega}$ skew-symmetric.
Hence, the space $\Pcal_r$ is isomorphic to the sum of the linear space of
translation and the linear space of skew-symmetric matrices, and its dimension is
therefore $d(d + 1)/2$. The space of non rigid body motion is the quotient of $\Pcal$
with respect to $\Pcal_r$, which we denote $\Pcal/\Pcal_r$. We introduce the
projection operator $\proj_c$ in $\Pcal$ defined as
\begin{equation}
  \label{eq:defprojc}
  \proj_c(\vect{l}) = \frac{1}{2}(\mat{L} + \mat{L}^\trans)(\vect{x} - \bar{\vect{x}}_E),
\end{equation}
for any $\vect{l}(\vect{x}) = \vect{u} + \mat{L}\vect{x}\in\Pcal$.  Here,
$\bar{\vect{x}}_E$ denotes the arithmetic average of the positions $\vect{x}_i$ of
all the nodes of the cell $E$, that is 
\begin{equation*}
  \bar{\vect{x}}_E = \frac{1}{n}\sum_{i=1}^{n}\vect{x}_i,
\end{equation*}
where $n$ corresponds to the number of nodes in $E$. We can check that $\proj_{c}$ is
a projection and $\proj_{c}(\vect{l}) = 0$ if and only if
$\vect{l}\in\Pcal_r$. Hence, the image of $\proj_c$ is in bijection with the space of
linear strain $\Pcal\quot\Pcal_r$ which we therefore identify to
$\Pcal_c=\proj_c(\Pcal)$. Note that $\Pcal_c$ can also be defined as
\begin{equation}
  \label{eq:altdefpcalc}
  \Pcal_c = \{\vect{l}\in\Pcal\ |\ \grad\vect{l}=\grad\vect{l}^\trans \text{ and } \vect{l}(\bar{\vect{x}}_E) = 0 \}.
\end{equation}
Then, we introduce the projection $\proj_r$ from $\Pcal$ to $\Pcal_r$ as
\begin{equation}
  \label{eq:defprojrpcal}
  \proj_r(\vect{l}) = \vect{l} - \proj_c(\vect{l})
\end{equation}
so that
\begin{equation*}
  \proj_r(\vect{l}) =  \vect{l}(\bar{\vect{x}}_E) + \frac12(\mat{L} - \mat{L}^\trans)(\vect{x} - \bar{\vect{x}}_E).
\end{equation*}
We can check that $\Pcal_r=\proj_r(\Pcal)$, $\proj_c\proj_r = \proj_r\proj_c = 0$,
$\proj_c + \proj_r = \id_{\text{\scalebox{0.7}{$\Pcal$}}}$.  The space $\Pcal_c$ is
isomorphic to the space of symmetric matrices, denoted $\Sym$. We consider the case
$d=3$ and use Kelvin's notation to represent $\Sym$ so that, for any
$\mat{a}\in\Sym$, its Kelvin representation in $\Real^6$, which we denote
$\hat{\vect{a}}\in\Real^6$, is given by
\begin{equation}
  \label{eq:defKelvnot}
  \hat{\vect{a}}^\trans = [a_{11}, a_{22}, a_{33}, \sqrt{2}a_{23}, \sqrt{2}a_{13}, \sqrt{2}a_{12}].
\end{equation}
The square root in the definition above has the advantage to lead to the following
correspondence between the scalar products in $\Sym$ and $\Real^6$,
\begin{equation}
  \label{eq:scalprodcor}
  \mat{a} : \mat{b} = \hat{\vect{a}}\cdot\hat{\vect{b}},
\end{equation}
for any $\mat{a},\mat{b}\in\Sym$. Note that the authors in \cite{gain2014} use Voigt
instead of Kelvin notations, which explains why the expressions given in the present
paper differ up to a coefficient to those in \cite{gain2014}. We define the symmetric
tensor $\hat{\mat{C}}\in\Real^{6\times 6}$ by the identity
\begin{equation}
  \label{eq:defhatC}
  \mat{a} : \mat{C}\mat{b} = \hat{\vect{a}}^\trans\hat{\mat{C}}\hat{\vect{b}},
\end{equation}
for all $\mat{a},\mat{b}\in\Sym$.  Then, we obtain that,
for any $\vect{l}, \vect{m}\in\Pcal$, we have
\begin{equation}
  \label{eq:energyinPcal}
  a_E(\vect{l}, \vect{m}) = \int_E\vect{\epsi}(\vect{l}):\mat{C}\vect{\epsi}(\vect{m})\dx = \abs{E}\widehat{\proj_c(\vect{l})}^\trans\hat{\mat{C}}\widehat{\proj_c(\vect{m})}
\end{equation}
The projection $\proj_c(\vect{l})$ belongs to $\Pcal_c$ and can therefore be written
as $\proj_c(\vect{l}) = \mat{Q}(\vect{x} - \hat{\vect{x}}_E)$ for $\mat{Q}\in\Sym$.
In \eqref{eq:energyinPcal}, we wrote $\widehat{\proj_c(\vect{l})}$ for the Kelvin
representation $\hat{\mat{Q}}$ of $\mat{Q}$, and similarly for
$\widehat{\proj_c(\vect{m})}$. For the space $\Pcal_r$, we use the mapping between
$\aSym$ and $\Real^3$ given by the cross-product operation. For any
$\mat{a}\in\aSym$, we can define a rotation vector $\hat{\vect{a}}\in\Real^3$ as
\begin{equation}
  \label{eq:defaasym}
  \hat{\vect{a}} = \sqrt{2}[-a_{2, 1}, a_{1, 2}, -a_{1, 2}]^\trans.
\end{equation}
Then, we have
\begin{equation}
  \label{eq:rotdep}
  \mat{a}\vect{x} = \frac 1{\sqrt{2}}\,\hat{\vect{a}}\times\vect{x}
\end{equation}
and
\begin{equation}
  \label{eq:equivscalrot}
  \mat{a}:\mat{b} = \hat{\vect{a}}\cdot\hat{\vect{b}}
\end{equation}
for any $\mat{a},\mat{b}\in\aSym$. The normalization using $\frac1{\sqrt 2}$ in
\eqref{eq:rotdep} is used in order to get an exact correspondence between the scalar
products in $\aSym$ and $\Real^3$, as in \eqref{eq:scalprodcor} in the case of
$\Sym$. The basis of $\Pcal_r$ is given by the canonical basis of $\Real^6$, the
first three components corresponding to the translation vector while the three last
components correspond to a rotation vector.

\subsection{The projection operator}
\label{subsec:projoperator}

The projection operator on the space of affine displacement, with respect to the
bilinear form $a$, plays an essential role in the formulation of a first order VEM
method. We follow the notation of \cite{beirao2013basic} and denote this projection
by $\projgrad$, even if the bilinear form we consider is not the $H^1$ semi-norm,
which is used as example in \cite{beirao2013basic}. For a given displacement function
$\vect{\nu}\in\Vcal_E$, the projection $\vect{p} = \projgrad(\vect{\nu})$ is defined
as the unique element $\vect{p}\in\Pcal$ which satisfies
\begin{subequations}
  \label{eq:defprojgrad}
  \begin{equation}
    \label{eq:defprojgrad1}
    a_E(\vect{p},\vect{q}) = a_E(\vect{\nu},\vect{q}),
  \end{equation}
  for all $\vect{q}\in\Pcal$ and such that 
  \begin{equation}
    \label{eq:defprojgrad2}
    \proj_r\vect{p} = \proj_r{\vect{\nu}},
  \end{equation}
\end{subequations}
which means that the projection of $\vect{p}$ and $\vect{\nu}$ on $\Pcal_r$
coincide. The condition \eqref{eq:defprojgrad2} is necessary to determine a unique
solution. Indeed, the bilinear form $a$ is degenerate as it is invariant with respect
to the space of rigid body motion and the condition \eqref{eq:defprojgrad2}
eliminates this underteminancy by imposing a rigid body motion for $\vect{p}$. It is
important to note that the definition of the projection operator is in fact
independent of $\Vcal_E$ and could be extended to $[H^1(E)]^3$. However, for
$\Vcal_E$, the projection can be computed exactly, using directly the degrees of
freedom and without further integration. Up to now, we have only introduced
$\proj_r$ and $\proj_c$ on $\Pcal$ and we will now extend these definitions to
$\Vcal_E$ so that the definition in \eqref{eq:defprojgrad2} makes actually sense. We
let the reader check that the new definitions of $\proj_c$ and $\proj_r$ on
$\Vcal_E$, when restricted to $\Pcal$, coincide with those introduced previously. We
define the projection $\proj_r: \Vcal_E \to \Pcal_r$ as
\begin{equation}
  \label{eq:defprojr}
  \proj_r(\vect{\nu}) = \bar{\vect{\nu}}_E + \frac12\brac{\nabla\vect{\nu} - \nabla\vect{\nu}^\trans} (\vect{x} - \bar{\vect{x}}_E),
\end{equation}
where the bracket denote the cell average, i.e.,
\begin{equation*}
  \brac{\mat{w}} = \frac{1}{\abs{E}}\int_E \mat{w}\dx\bigand\bar{\vect{\nu}}= \frac{1}{n}\sum_E \vect{\nu}_i.
\end{equation*}
We define the projection $\proj_c:\Vcal_E\to\Pcal_c$ as
\begin{equation*}
  \proj_c(\vect{\nu}) = \frac12\brac{\grad\vect{\nu} + \grad\vect{\nu}^\trans}(\vect{x} - \bar{\vect{x}}_E).
\end{equation*}
In a moment, we are going to check that both projections can be computed directly
from the degree of freedoms. First, we use these definitions to compute $\projgrad$.
We start by considering a solution $\vect{p}=\projgrad(\vect{\nu})$ to
\eqref{eq:defprojgrad} and show that \eqref{eq:defprojgrad1} yields  
\begin{equation}
  \label{eq:scalprodred}
  a_E(\vect{p}_c,\vect{q}_c) = a_E(\vect{\nu},\vect{q}_c),
\end{equation}
for any $\vect{q}\in\Pol_1(E)$, where $\vect{p}_c = \proj_c(\vect{p})$ and
$\vect{q}_c = \proj_c(\vect{q})$. The symmetric gradient is zero for any element in
$\Pcal_r$, that is $\epsi\circ\proj_r = 0$. Hence, $a_E(\proj_r(\vect{p}),
\vect{\nu}) = 0$ for any $\vect{p}\in\Pol_1(E)$ and $\vect{\nu}$. It implies that
$a_E(\vect{p},\vect{q}) = a_E(\vect{p}_c,\vect{q}_c) $ and $a_E(\vect{\nu},\vect{q})
= a_E(\vect{\nu},\vect{q}_c)$ so that Equation \eqref{eq:defprojgrad1} indeed implies
\eqref{eq:scalprodred}. Let us now determine the element $\vect{p}$ that satisfies
\eqref{eq:defprojgrad} for a given $\vect{\nu}\in\Vcal_E$. The coercivity of the the
form $a_E$ on $\Pcal_c$ follows from the definition of $\Pcal_c$ and the coercivity
of the tensor $\mat{C}$, see \eqref{eq:elliptcond}. Therefore, there exists a unique
solution $\vect{p}_c\in\Pcal_c$ such that \eqref{eq:scalprodred} holds for all
$\vect{q}_c\in\Pcal_c$. For any $\vect{q}_c\in\Pcal_c$, we have
\begin{equation}
  \label{eq:aenuqc}
  a_E(\vect{p}_c,\vect{q}_c) = \int_E \grad\vect{p}_c:\mat{C}\grad\vect{q}_c\dx = \abs{E}\grad\vect{p}_c:\mat{C}\grad\vect{q}_c
\end{equation}
and
\begin{equation}
  \label{eq:aenuqc}
  a_E(\vect{\nu},\vect{q}_c) = \int_E \frac12(\grad\vect{\nu} + \grad\vect{\nu}^\trans):\mat{C}\grad{\vect{q}}_c\dx = \left(\frac12\int_E(\grad\vect{\nu} + \grad\vect{\nu}^\trans)\dx\right):\mat{C}\grad{\vect{q}}_c.
\end{equation}
Hence, $\grad\vect{p}_c = \frac12\int_E(\grad\vect{\nu} + \grad\vect{\nu}^\trans)\dx$
which implies that $\vect{p}_c$ is uniquely defined as $\vect{p}_c =
\proj_c(\vect{\nu})$. We can conclude that $\vect{p}$ defined as
\begin{equation}
  \vect{p} = \vect{p}_c + \proj_r(\vect{\nu})
\end{equation}
is the unique solution to \eqref{eq:defprojgrad}. Indeed, $\proj_c(\vect{p}) =
\vect{p}_c$ and $\proj_r(\vect{p}) = \proj_r(\vect{\nu})$ are both uniquely defined
by \eqref{eq:scalprodred} and \eqref{eq:defprojgrad2}.
 
Let us now give more details on the assembly. To do so, we consider a basis function
$\vect{\nu}^i\in\Vcal_E$ for which the only non-zero displacement can only occur at
the node $i$, that is $\vect{\nu}^i(\vect{x}_j) = 0$ if $i\neq j$. Such function can
be written as
\begin{equation*}
  \vect{\nu}^i(\vect{x}) = \textstyle\sum_{j=1}^d\nu_j\phi_i(\vect{x})\vect{e}_j,
\end{equation*}
where $\{\vect{e}\}_{j=1}^{d}$ is the basis for Cartesian coordinates. We have
\begin{equation*}
  \brac{\grad\vect{\nu}^i} = \textstyle\sum_{j=1}^d\nu_j^i\vect{e}_j\brac{\grad\phi_i}^\trans,
\end{equation*}
and we have to compute $\vect{q}^i = \brac{\grad\phi_i}$. The expression above
simplifies to
\begin{equation*}
  \brac{\grad\vect{\nu}^i} = \vect{\nu}^i\vect{q}^{i\trans}.
\end{equation*}
For $\vect{q}^i$, using Stokes' theorem, we have
\begin{equation}
  \label{eq:defqi}
  \vect{q}^i = \int_E\grad\phi_i\dx= \int_{\partial E} \phi_i\vect{n} \dx = \sum_{f\in F(E)} (\int_{f}\phi_i \dx)\,\vect{n}_f,
\end{equation}
where $F(E)$ denotes the set of faces that belong to $E$. The integral in
\eqref{eq:defqi} can be computed exactly. For the 3D, we use a virtual space such
that the first two moments of the virtual basis elements coincide with those of their
projection, see \cite{ahmad2013equivalent}. The integral is zero if the node $i$ does
not belong to the face $f$ and, otherwise,
\begin{equation}
  \label{eq:explfaceint}
  \int_{f}\phi_i \dx =
  \begin{cases}
    \frac{|f|}{m} +\frac{1}{2}(\vect{n}_{e,i^-}+\vect{n}_{e,i^+})\cdot (\vect{x}^f-\bar{\vect{x}}^f) & \text{ in 3D},\\
    \frac{|f|}{2} & \text{ in 2D},
  \end{cases}
\end{equation}
where $\vect{x}^f$ is the centroid of the face $f$ and $\bar{\vect{x}}^f$ the
arithmetic average of the node coordinates, i.e. $\bar{\vect{x}}^f =
\frac1{m}\sum_{j=1}^m\vect{x}_j^f$. We denote by $\mat{W}_c^i\in\Real^{6\times 3}$
the matrix representation of $\proj_c$ written in the basis of displacement for the
node $i$ (that is $\Real^3$) and the basis of $\Pcal_c$ (that is $\Real^6$, using the
Kelvin notation). For $l,m = \{1, 2, 3\}$, we have
\begin{equation*}
  \frac12\brac{\grad\mat{\nu}^i + \grad\mat{\nu}^{i\trans}}_{l,m} = \frac12(\nu_l^iq_m^i + \nu_m^iq_l^i)
\end{equation*}
so that
\begin{equation*}
  (\mat{W}_c^{i})^\trans =
  \begin{pmatrix}
    q_1^i&0&0&0&\smallinvroottwo q_3^i&\smallinvroottwo q_2^i\\
    0&q_2^i&0&\smallinvroottwo q_3^i&0&\smallinvroottwo q_1^i\\
    0&0&q_3^i&\smallinvroottwo q_2^i&\smallinvroottwo q_1^i&0
  \end{pmatrix}.
\end{equation*}
We have
\begin{equation*}
  \frac12\brac{\grad\mat{\nu}^i - \grad\mat{\nu}^{i\trans}} = \frac12(\vect{\nu}^i\vect{q}^{i\trans} - \vect{q}^i\vect{\nu}^{i\trans}).
\end{equation*}
Using the general identity $(\vect{q}^i\times\vect{\nu}^i)\times\vect{x} =
(\vect{q}^i\cdot\vect{x})\vect{\nu}^i - (\vect{\nu}^i\cdot\vect{x})\vect{q}^i$, we
get that the $\Real^3$ representation of the matrix in $\aSym$ above is given by
$\frac1{\sqrt 2}\vect{q}^i\times\vect{\nu}^i$.  Hence, the matrix
$\mat{W}_r^i\in\Real^{6\times 3}$ that represents $\proj_r$ written in the basis of
displacement for the node $i$ (that is $\Real^3$) and the basis of $\Pcal_r$ (that is
$\Real^6$ for the translation and the rotation vector) is given by
\begin{equation*}
  (\mat{W}_r^{i})^\trans =
  \begin{pmatrix}
    \frac1n&0&0&0&\smallmininvroottwo q_3^i&\smallinvroottwo q_2^i\\
    0&\frac1n&0&\smallinvroottwo q_3^i&0&\smallmininvroottwo q_1^i\\
    0&0&\frac1n&\smallmininvroottwo q_2^i&\smallinvroottwo q_1^i&0
  \end{pmatrix}.
\end{equation*}
The matrices $\mat{W}_c$ from the space of all the degrees of freedom (that is
$\Real^{3n}$) to $\Pcal_c$ is obtained by concatenating $\mat{W}_c^i$ and similarly
for $\mat{W}_r$. To obtain, from $\mat{W}_c$ and $\mat{W}_r$, the matrix
representations of $\proj_c$ and $\proj_r$ in terms only of the degrees of freedom,
we have to find the decomposition of $\Pcal_c$ and $\Pcal_r$ in terms of the degrees
of freedom. To do so, we introduce the vectors $\vect{r}^i$, for $i = \{1,\ldots,
n\}$ as
\begin{equation*}
  \vect{r}^i = \vect{x}_i - \bar{\vect{x}}_E.
\end{equation*}
We define $\mat{N}_c^i,\mat{N}_c^i\in\Real^{3\times 6}$ as
\begin{equation}
  \label{eq:defNicNir}
  \mat{N}_c^{i} =
  \begin{pmatrix}
    r_1^i&0&0&0&\smallinvroottwo r_3^i&\smallinvroottwo r_2^i\\
    0&r_2^i&0&\smallinvroottwo r_3^i&0&\smallinvroottwo r_1^i\\
    0&0&r_3^i&\smallinvroottwo r_2^i&\smallinvroottwo r_1^i&0
  \end{pmatrix}\text{ and }\mat{N}_r^{i} =
  \begin{pmatrix}
    1&0&0&0&\smallmininvroottwo r_3^i&\smallinvroottwo r_2^i\\
    0&1&0&\smallinvroottwo r_3^i&0&\smallmininvroottwo r_1^i\\
    0&0&1&\smallmininvroottwo r_2^i&\smallinvroottwo r_1^i&0
  \end{pmatrix}.
\end{equation}
Then, the matrices $\mat{N}_c,\mat{N}_r\in\Real^{3n\times 6}$ are obtained by
concatenating $\mat{N}_c^i,\mat{N}_r^i$, respectively.

The projections can be then written in terms as a mapping from degrees of freedom to
degrees of freedom,
\begin{equation}
  \dproj_r = \mat{N}_r\mat{W}_r  \quad\text{ and }\quad\dproj_c = \mat{N}_c\mat{W}_c
\end{equation}
and the projection on affine displacement is given by $\dproj = \dproj_c +
\dproj_r$. For any $\vect{\nu}, \vect{\eta}\in\Vcal_E$, we have that
\begin{equation*}
  a_{E}(\projgrad\vect{\nu}, \projgrad\vect{\eta}) = \abs{E}\widehat{\proj_c(\vect{\nu})}^\trans\hat{\mat{C}}\widehat{\proj_c(\vect{\eta})} = \abs{E}\vect{\nu}^\trans\vect{W_c}^\trans\hat{\mat{C}}\mat{W}_c\vect{\eta},
\end{equation*}
where in the last term, slightly abusing the notations, we denote by
$\vect{\nu},\vect{\eta}\in\Real^{3n}$ the vector composed of the degrees of freedoms
of $\vect{\nu},\vect{\eta}\in\Vcal_E$. Using the same convention, we can write the
bilinear form as
\begin{equation}
  \label{eq:ahexpr}
  a_{h,E}(\vect{\nu}, \vect{\eta}) =  \vect{\nu}^\trans\left(\abs{E}\vect{W_c}^\trans\hat{\mat{C}}\mat{W}_c+ (\id-\dproj)^T\mat{S}(\id-\dproj)\right)\vect{\eta},
\end{equation}
where $\mat{S}\in\Real^{n\times n}$ is a \textit{stabilization} term. For the VEM
method to be well-defined, the matrix $\mat{S}$ must be chosen such that it is
positive, symmetric and definite on the kernel of $\mat{P}$. We note that the
decomposition of the energy in two orthogonal parts, the linear part which ensures
consistency and the higher order part which are handled so that stability is
preserved, is analog to the decomposition used in \cite{Bergan1984}, even if it was
introduced there to add some freedom in the choice of the basis functions.

\section{Implementation of the load term}
\label{subsec:loadterm}

The load term can be calculated in several different ways which are equivalent up to
the order of accuracy of the methods. We have investigated the three following
alternatives,
\begin{enumerate}
\item Computation using the projection operator $\projgrad$,
\item Integration using nodal quadrature,
\item Computation based on a discrete gradient operator.
\end{enumerate}
Alternative 1 is the choice that naturally follows from the VEM approach and which is
proposed in \cite{beirao2013basic}. Alternative 2 was argued to be simpler and with
similar accuracy in \cite{gain2014}. Alternative 3 is possible when the force is
equal to the gradient of a potential. This last alternative actually came to the mind
of the authors when they considered the poro-elasticity equation, where the
divergence operator naturally arises. As we will see below, the discrete gradient is
in fact derived from the discrete divergence operator by duality. The two first
alternatives give similar results. We show that the last one has significantly less
errors than the others for elongated grid cells.

\subsection{Standard assembly of the load term (alternative 1 and 2)}

For a given force $\vect{f}$, we consider the work done by the force for a given
displacement field $\vect{u}$,
\begin{equation}
  \label{eq:linformworkdef}
  \int_\Omega\vect{f}\cdot\vect{u}\dx.
\end{equation}
This expression defines a linear form on the space of displacement. We denote by
$\Vcal$ the global discrete function space of displacement, which is constructed by
taking the product of the $\Vcal_E$ for all the cells $E$ of the grid and requiring
continuity at the cell boundaries and correspond to the nodal displacement in terms
of the degrees of freedom. We want to find a discrete linear form on $\Vcal$ that
approximates \eqref{eq:linformworkdef}. We can equip $\Vcal$ with the standard scalar
product in $\Real^{n_N}$, that is
$\sum_{\eta}\vect{u}_\eta\cdot\vect{\nu}_\eta$. Here, $n_N$ denotes the total number
of nodes. Any linear form on $\Vcal$ can be represented by an element in $\Vcal$,
using this scalar product. Hence, we end up looking for an element
$\widehat{\vect{f}}\in\Vcal$ such that
\begin{equation}
  \label{eq:nodalforce}
  \int_\Omega\vect{f}\cdot\vect{u}\dx\approx \sum_{\eta} \widehat{\vect{f}}_\eta\cdot\vect{u}_\eta.
\end{equation}
The vector $\widehat{\vect{f}}\in\Real^{N_n}$ can be interpreted as a vector of
\textit{nodal forces}. We present several expressions for $\widehat{\vect{f}}$
corresponding to the three alternatives presented previously. First, we can use
weights which are obtained using a first-order quadrature. For a node $\eta$, let us
denote by $E(\eta)$ the set of cells to which the node $\eta$ belongs. Using
quadrature rules to integrate $f$ on each cell, we obtain
\begin{equation}
  \widehat{\vect{f}}_\eta = \left(\sum_{i\in E(\eta)}w_i^\eta\right)\vect{f}(\eta),
\end{equation}
see \cite{gain2014} for the definitions of the weights $w_{i}^\eta$. This corresponds
to alternative 2. For alternative 1, let $\vect{u}^\eta\in\Vcal$ be a displacement
for which the only non-zero degrees of freedom are those corresponding to the node
$\eta$. Then, we have
\begin{align}
  \notag
  \int_\Omega\vect{f}\cdot\vect{u}^{\eta}\dx &= \sum_{i\in E(\eta)}\int_{E_i}\vect{f}\cdot\vect{u}^{\eta}\dx\\
  \notag
  &\approx \sum_{i\in E(\eta)}\int_{E_i}\proj_i^0(\vect{f})\cdot\vect{u}^{\eta}\dx\\
  \notag
  &= \sum_{i\in E(\eta)}\proj_i^0(\vect{f})\cdot\int_{E_i}\vect{u}^{\eta}\dx\\
  \label{eq:exactf1}
  &= \sum_{i\in E(\eta)}\proj_i^0(\vect{f})\cdot\int_{E_i}\projgradi(\vect{u}^{\eta})\dx.
\end{align}
Here, $\proj_i^0$ denotes the $L^2$ projection to the space of constant functions
(polynomials of degree zero) in the element $E_i$.  To obtain the last integral, we
use that fact that the virtual basis functions can be chosen such that the zero and
first moment of a function $\vect{\nu}$ in $\Vcal_E$ coincide with those of its
projection $\projgradi\vect{\nu}$, that is
\begin{equation*}
  \int_{E}\vect{p}\cdot\vect{\nu}\dx = \int_{E}\vect{p}\cdot\projgradi(\vect{\nu})\dx,
\end{equation*}
for any $\vect{p}\in\Pol_1$ and $\vect{\nu}\in\Vcal_E$. See
\cite{ahmad2013equivalent} for more details. The choice of such basis implies that,
for an element $E$, the modes that belong to $\ker\projgrad$, typically higher
nonlinear modes, will not be excited directly by the force. From \eqref{eq:exactf1},
we infer that $\widehat{\vect{f}}_\eta$ is defined as a linear combination of cell
averages of $\vect{f}$,
\begin{equation}
  \widehat{\vect{f}}_{\eta} = \sum_{i\in E(\eta)} m_{i}^{\eta} \proj_i^0(\vect{f})
\end{equation}
where $m_i^\eta$ are the weights given
\begin{equation*}
  m_{i}^{\eta} =\vect{e}_k\cdot\int_{E_i}\projgradi(\vect{\phi}_{k}^\eta)\dx.
\end{equation*}
Note that the expression on the right above do not depend on $k$, as the same basis
function is used in all directions.

\subsection{The discrete gradient approach}

Let us now turn to alternative 3. We assume that the force can be written as the
gradient of a potential, $\vect{f} = \nabla \psi$. We have, for a node $\eta$ and a
dimension $k\in\{1,2,3\}$,
\begin{equation}
  \label{eq:workasdiv}
  \int_{\Omega} \vect{f}\cdot\vect{u} \dx = \int_{\Omega} \grad \psi \cdot
  \vect{u} \dx = - \int_{\Omega} \psi \dive \vect{u} \dx + \int_{\partial \Omega}
  \psi
  \vect{u}\dx
  = - \int_{\Omega} \psi \dive \vect{u} \dx.
\end{equation}
The boundary integral vanishes because we assume Dirichlet boundary condition,
$\vect{u} = 0$ for $\vect{x}\in\partial\Omega$. In the VEM space, there exists a
natural discretization of the divergence operator as an operator from $\Vcal$ to
cell-wise constant functions, denoted $\Tcal$, which is isomorphic to $\Real^{N_c}$,
where $N_c$ denotes the number of cells. Indeed, for any discretized potential
$\hat\psi\in\Tcal$ and $\vect{\nu}\in\Vcal$, we have
\begin{equation}
  \label{eq:globddiv}
  \int_\Omega \hat\psi\dive\vect{\nu} = \sum_{E}\int_E\hat\psi_E\dive\vect{\nu}\dx = \sum_{E}\sum_{f\in F(E)}\hat\psi_E\int_f\vect{\nu}\cdot\vect{n}\dx,
\end{equation}
where $F(E)$ as before denotes the set of faces that belong to $E$. The last integral
can be computed exactly as shown in \eqref{eq:explfaceint}.  Then, using partial
integration, we get
\begin{equation}
  \label{eq:ddivdef}
  \int_\Omega\hat\psi\dive\vect{\phi}_{\eta}^{k}\dx = \sum_{j\in E(\eta)}\sum_{f_{j,l}\in E_j\cap E_l}\hat\psi_{E_j}(\vect{e}_k\cdot\vect{n}_{j.l})\int_{f_{j,l}}\phi_{\eta}\dx,
\end{equation}
with the convention that we only get contribution in the integral when the face
$f_{j,l}$ exists, that is when $E_j$ and $E_l$ share a common face. Note that by
definition of the exterior normal, we have $\vect{n}_{j.l} = -\vect{n}_{j.l}$. We use
\eqref{eq:explfaceint} to compute the integral and therefore the divergence operator
$\ddiv:\Vcal\to\Tcal$ is defined and can be computed exactly in the sense that 
\begin{equation*}
  \ddiv(\vect{\nu}) = \proj^0\dive\vect{\nu}
\end{equation*}
for any $\vect{\nu}\in\Vcal$, where $\proj^0$ denotes the $L^2$ projection to $\Tcal$. The transpose of the discrete divergence operator will
give us a discrete approximation of the gradient. We can obtain an expression of the
discrete gradient by reverting the order of the sum in \eqref{eq:ddivdef}. Let us
denote by $F(\eta)$ the set of faces to which the node $\eta$ belongs and, for a face
$f_k$, we denote the neighboring cells of $f_k$ by $E_{k}^+$ and $E_{k}^-$. From
\eqref{eq:ddivdef}, we can rewrite
\begin{equation*}
  \int_\Omega\hat\psi\dive\vect{\phi}_{\eta}^k\dx = -\sum_{f\in F(\eta)}(\hat\psi_{E_f^+} - \hat\psi_{E_f^-})(\vect{e}_k\cdot\vect{n}_{f})\int_{f}\phi_{\eta}\dx.
\end{equation*}
where the normal $\vect{n}_f$ of the face $f$ is directed from $E_f^-$ to
$E_f^+$. This convention implies that $\vect{n}_{jl}=-\vect{n}_{lj}=\vect{n}_f$ if
$E_j=E_f^-$ and $E_l=E_f^+$. Hence, the discrete gradient operator $\dgrad$ is the
mapping from scalar cell values to vector node value given by
\begin{equation}
  \label{eq:defdgrad}
  [\dgrad(\hat\psi)]_{\eta, k} = \sum_{f\in F(\eta)}(\hat\psi_{E_f^+} - \hat\psi_{E_f^-})(\vect{e}_k\cdot\vect{n}_{f})\int_{f}\phi_{\eta}\dx.
\end{equation}
Hence, gathering \eqref{eq:workasdiv}, \eqref{eq:globddiv} and \eqref{eq:defdgrad},
in this formulation, we obtain the following expression for $\widehat{\vect{f}}$, as
the discrete gradient of the discretized potential, that is
\begin{equation*}
  \widehat{\vect{f}} = \dgrad(\hat\psi).
\end{equation*}
In \eqref{eq:defdgrad}, the expression only depends on differences of the potential,
which can be estimated locally without knowledge of the global potential, i.e.
\begin{equation}
  \hat\psi_{E_f^+} - \hat\psi_{E_f^-}= \hat{\vect{f}}_f\cdot \vect{dr}_f,
  \label{eq:local_faceforce}
\end{equation}
where $\hat{\vect{f}}_f$ is an approximation of the force on the face $f$ and
$\vect{dr}_f$ is the vector joining the centroids of ${E_f^-}$ and ${E_f^+}$. In
practice, it means that the method can be applied \textit{even} if the force is not
derived from a potential, as we can see that the potential $\psi$ does not have to be
computed. Note that, in the numerical tests that follow, we have not tested this
case.

\subsection{Interpretation of the discrete gradient approach using singular load
  term functions}

When we consider a cell-valued potential $\psi$, the corresponding force
$\vect{f}=\grad\psi$ can be defined as a singular function with support on the cell
faces.  Let us define this class of function, which we will refer to as
\textit{2D-Dirac} functions. Given an internal 2D surface $S$ in $\Omega$ (or 1D line
in 2D), we define the \textit{constant 2D-Dirac} function $\delta_S(x)$ as the
distribution given by
\begin{equation*}
  <\delta_S,\phi> = \int_S\phi(x)\,dx,
\end{equation*}
for all $\phi\in C^\infty(\Omega)$. The 2D-Dirac $\delta_S$ is a measure which
coincides with the Hausdorff measure on the $d-1$ dimensional set $S$. Then, we can
also define 2D-Dirac function $h(x)\delta_S(x)$, for any $h\in L^1(S)$ as
$<h\delta_S,\phi> = \int_Sh(x)\phi(x)\,dx,$, for any $\phi\in C^\infty(\Omega)$. If
the surface $S$ is Lipschitz, then a continuous trace operator from $H^1(\Omega)$ to
$H^{\frac12}(S)$ can be defined, see for example \cite{ding1996proof}. Therefore, at
least if $h\in L^2(S)$, we have that $h(x)\delta_S(x)\in H^{-1}(\Omega)$. Indeed, we
have, for any $\phi\in H^1(\Omega)$,
\begin{align*}
  <h\delta_S,\phi> = \int_S h(x)\phi(x)\,dx &\leq \norm{h}_{L^2(S)}\norm{\phi}_{L^2(S)}\\
                                            &\leq C_1 \norm{h}_{L^2(S)}\norm{\phi}_{H^{\frac12}(S)}\\
                                            &\leq C_1C_2 \norm{h}_{L^2(S)}\norm{\phi}_{H^1(\Omega)},
\end{align*}
for two constants $C_1$ and $C_2$. From this observation, we can infer that the
original system of equation \eqref{eq:lin_elast_cont} is \textit{well-posed} for 2D
Dirac vector functions $\vect{f}$.

Let us now consider a surface $S$ that splits the domain $\Omega$ in two sub-domains,
namely $\Omega_-$ and $\Omega_+$, and a potential $\psi$ which is piecewise constant
and takes the values $\psi_{\pm}$ in $\Omega_{\pm}$. The gradient of $\psi$ in the
sense of distribution is defined as
\begin{equation}
  \label{eq:defweakgrad}
  <\grad\psi,\phi> = - \int \psi\grad\phi \,dx 
\end{equation}
Let us consider $\phi$ with compact support in $\Omega$ so that, after using
integration by part we obtain,
\begin{equation}
  \label{eq:weakgradcalc}
  <\grad\psi,\phi> = - (\psi_-\int_{\Omega_-} \grad\phi \,dx + \psi_+\int_{\Omega_+} \grad\phi \,dx) =
  =\int_S((\psi_+ - \psi_-)\vect{n})\phi\,dx,
\end{equation}
where $\vect{n}(x)$ denotes the normal to $S$ at $x\in S$ pointing from $\Omega_-$ to
$\Omega_+$. From the definition \eqref{eq:defweakgrad} and \eqref{eq:weakgradcalc},
we get that the gradient of $\hat{\psi}$ is a 2D Dirac vector function given by
\begin{equation}
  \label{eq:weakgradpot}
  \grad\psi = [\psi_+ - \psi_-]\vect{n}(x)\delta_S.
\end{equation}

Let us now consider again a cell-wise constant potential function $\hat\psi$ defined
on a mesh. Using the same notation as in the previous section, we infer from
\eqref{eq:weakgradpot} that the gradient $\mathring{\vect{f}}$ of $\hat\psi$ in the
sense of distribution is given by
\begin{equation}
  \label{eq:defder2ddirac}
  \mathring{\vect{f}} = \grad\hat\psi = \sum_{f\in \Fint}(\hat\psi_{E_f^+} - \hat\psi_{E_f^-})\vect{n}_f\delta_f,
\end{equation}
where $\Fint$ denotes the set of internal faces. Note that for the basis function
$\vect{\phi}_{\eta}^{k}$ as defined in \eqref{eq:defbasis}, we get
\begin{equation}
  \label{eq:exp2dirass}
  \int_\Omega \mathring{\vect{f}}\cdot\vect{\phi}_{\eta}^{k} = \sum_{f\in F_i}(\hat\psi_{E_f^+} - \hat\psi_{E_f^-})(\vect{e}_k\cdot\vect{n}_f)\int_f\phi_\eta(x)\,dx
\end{equation}
and we recover expression \eqref{eq:defdgrad}. Hence, the discrete gradient approach
can be interpreted in the following way. First, we approximate the volumetric load
term $\vect{f}$ by a 2D Dirac function $\mathring{\vect{f}}$ with support on the cell
faces and which is constant on each face, that is $\mathring{\vect{f}}$ has the form
\begin{equation}
  \label{eq:exprconstf}
  \mathring{\vect{f}}(x) = \sum_{f\in F_i} \mathring{\vect{f}}_f\delta_f(x)
\end{equation}
where $\mathring{\vect{f}}_f$ is a constant vector, for each face $f$. In the case
the force $\vect{f}$ is derived from a potential, we can use the expression
\eqref{eq:defder2ddirac} to carry on this approximation. Otherwise, we propose to use
expression \eqref{eq:local_faceforce} and consider
\begin{equation*}
  \mathring{\vect{f}}_f =  (\hat{\vect{f}}_f \cdot \vect{dr}_f)\,\vect{n}_f.
\end{equation*}
Once $\mathring{\vect{f}}$ is computed, we use the VEM method to solve the problem defined
as
\begin{equation*}
  \dive \mat\sigma = \mathring{\vect{f}}.
\end{equation*}
Then, the assembly of the load term can be done \textit{exactly}, as we can see from
\eqref{eq:exp2dirass} in the case of a potential and otherwise
\begin{equation}
  \label{eq:exprfsurf}
  \int_\Omega \mathring{\vect{f}}\cdot\vect{\phi}_{\eta}^{k} = \sum_{f\in F_i}\mathring{\vect{f}}_f\cdot\vect{e}_k\int_f\phi_\eta(x)\,dx.
\end{equation}
in the case where \eqref{eq:exprconstf} is used. Note that the integrals in
\eqref{eq:exp2dirass} and \eqref{eq:exprconstf} can be computed exactly we use the
virtual basis proposed in \cite{ahmad2013equivalent}.

\section{Stability with respect to aspect ratio}
\label{subsec:stabaspectrat}

Let us now discuss the choice of the stabilization matrix $\mat{S}$ in
\eqref{eq:ahexpr}. In \cite{da2013virtual}, the authors propose
\begin{equation*}
  \mat{S} = \id,
\end{equation*}
which is the simplest choice. In \cite{gain2014}, the authors look at several cell
shapes and recommend the stabilization term given by
\begin{equation}
  \label{eq:stabgain}
  \mat{S} = \alpha \id
\end{equation}
where the constant $\alpha$ is chosen as
\begin{equation}
  \label{eq:defalphagain}
  \alpha_G = \frac{\abs{E}\trace(\hat{\mat{C}})}{\trace(\mat{N}_c^\trans\mat{N}_c)},
\end{equation}
as it gives an overall satisfactory approximation of the higher order nonlinear
modes. This constant is stable with respect to isotropic scaling but it is not stable
with respect to the aspect ratio.

\subsection{Instability of $\alpha_G$ with respect to aspect ratio}

To demonstrate that, we consider a rectangular element in 2D given by $[-h_1,
h_1]\times[-h_2, h_2]$. In this case, an explicit definition of the virtual element
space is available, as it is spanned by the four following functions
\begin{equation}
  \label{eq:basiselt2d}
  \vphi_1^{l}(\vect{x}) = 1,\quad \vphi_2^{l}(\vect{x}) =  \frac{x_1}{h_1},\quad \vphi_3^{l}(\vect{x}) =  \frac{x_2}{h_2},\quad \vphi(\vect{x}) =  \frac{x_1x_2}{h_1h_2},
\end{equation}
in each Cartesian direction. They coincide in this case to the standard finite
elements for quadrilaterals. The functions $\vphi_j^l(\vect{x})\vect{e}_i$ for
$j=1,2,3$ and $i=1,2$ provides a basis for the affine space $\Pcal$. Let
$\vect{\vphi}_i(\vect{x}) = \vphi(\vect{x})\vect{e}_i$. We check directly, using the
symmetry of the domain, that
\begin{equation*}
  \projgrad(\vect{\vphi}_i) = 0.
\end{equation*}
Hence, for each basis functions in \eqref{eq:basiselt2d}, we have that the zero and
first order moments correspond to those of their projection so that, indeed, they
form a basis of $\Vcal$. Moreover $\{\vect{\vphi}_i\}_{i=1,2}$ constitutes a basis
for $\ker\projgrad$. In this two-dimensional case, the matrix $\mat{N}_c$ is given by
\begin{equation*}
  \mat{N}_c^{i} =
  \begin{pmatrix}
    h_1&0&\smallinvroottwo h_2\\
    0&h_2&\smallinvroottwo h_1
  \end{pmatrix}
\end{equation*}
We collect the contributions of the four nodes of the cell and obtain the matrix
$\mat{N}_c$ given by
\begin{equation}
  \label{eq:fulnc}
  \mat{N}_c^\trans =
  \begin{pmatrix}
    h_1                 &0                    & -h_1                & 0                     &     -h_1                 &0                    & h_1                & 0                      \\
    0                   &h_2                  & 0                   & h_2                   &     0                   &-h_2                  & 0                   & -h_2                    \\
    \smallinvroottwo h_2& \smallinvroottwo h_1& \smallinvroottwo h_2& -\smallinvroottwo h_1  &    -\smallinvroottwo h_2& -\smallinvroottwo h_1& -\smallinvroottwo h_2& \smallinvroottwo h_1
  \end{pmatrix}
\end{equation}
which yields
\begin{equation}
  \label{eq:ntn}
  \mat{N}_c^\trans\mat{N}_c =
  \begin{pmatrix} 4 h_1^2 & 0 &
    0 \\ 0 & 4 h_2^2 & 0 \\ 0 & 0 & 2(h_1^2 + h_2^2)
  \end{pmatrix}
\end{equation}
so that $\trace(\mat{N}_c^{\trans}\mat{N}_c) = 6(h_1^2 + h_2^2)$. Hence, the scaling
ratio $\alpha_G$ is given by
\begin{equation}
  \label{eq:alphaforrec}
  \alpha_G = \frac{4h_1h_2\trace(\hat{\mat{C}})}{6(h_1^2 + h_2^2)} = \frac{2}{3}\frac{\trace(\hat{\mat{C}})}{(\epsi + \epsi^{-1})},
\end{equation}
where $\epsi = \frac{h_1}{h_2}$ denotes the aspect ratio. Let us now compute how this
weight in the stabilization term compares with the actual energy for the functions
that belong to $\ker\projgrad$. To do so, we consider an isotropic material where the
stress is given as
\begin{equation}
  \label{eq:hook}
  \mat{\sigma} = \lambda\trace(\mat{\epsi}) + 2\mu\mat{\epsi},
\end{equation}
which implies
\begin{equation}
  \label{eq:aforiso}
  a(\vect{u}, \bar{\vect{u}}) = \int_\Omega (\lambda\trace(\mat{\epsi})\trace(\bar{\mat{\epsi}}) + 2\mu\mat{\epsi}:\bar{\mat{\epsi}} )\dx
\end{equation}
For $\vect{\vphi}_i$ we denote by $\mat{\epsi}_i$, the corresponding strain, which is
given by
\begin{equation}
  \label{eq:defepsii}
  \mat{\epsi}_i =   \frac12(\vect{e}_i\grad\phi^\trans + \grad\phi\vect{e}_i^\trans).
\end{equation}
We get
\begin{equation*}
  \mat{\epsi}_i : \mat{\epsi}_j = \frac12\left(\delta_{i,j}\abs{\grad\vphi}^2 + \fracpar{\vphi}{x_i}\fracpar{\vphi}{x_j}\right),
\end{equation*}
where $\delta_{i,j}=1$ if $i=j$ and zero otherwise. Hence, using the symmetry of the
domain, we get
\begin{equation*}
  \int_E \mat{\epsi}_i : \mat{\epsi}_j \dx = \delta_{i,j}\int_E\abs{\grad\vphi}^2\dx.
\end{equation*}
We have $\trace(\mat{\epsi}_i) = \fracpar{\vphi}{x_i}$. Hence, using the symmetry of
the domain we get
\begin{equation*}
  \int_E \trace(\mat{\epsi}_i)\trace(\mat{\epsi}_j)\dx = \delta_{i,j}\int_E\abs{\fracpar{\vphi}{x_i}}^2\dx.
\end{equation*}
Finally, the restriction of the bilinear form $a$ to $\ker\projgrad$ takes the form
\begin{equation*}
  a(\vect{\vphi}_i, \vect{\vphi}_j) =
  \begin{pmatrix}
    \int_E(\lambda\abs{\fracpar{\vphi}{x_1}}^2 + 2\mu\abs{\grad\vphi}^2)\dx&0\\
    0&\int_E(\lambda\abs{\fracpar{\vphi}{x_2}}^2 + 2\mu\abs{\grad\vphi}^2)\dx
  \end{pmatrix}
\end{equation*}
The integrals above can be computed exactly and we have
\begin{equation*}
  \int_E\abs{\fracpar{\vphi}{x_1}}^2\dx = \frac{4}{3}\epsi^{-1},\quad \int_E\abs{\fracpar{\vphi}{x_2}}^2\dx = \frac{4}{3}\epsi,
\end{equation*}
Hence, 
\begin{equation*}
  a\left(\vect{\vphi}_i, \vect{\vphi}_j\right) =
  \begin{pmatrix}
    \frac{4}{3}\lambda\epsi^{-1} + \frac{8}{3}\mu(\epsi + \epsi^{-1})&0\\
    0&\frac{4}{3}\lambda\epsi + \frac{8}{3}\mu(\epsi + \epsi^{-1})
  \end{pmatrix}.
\end{equation*}
We denote by $\alpha_1$ and $\alpha_2$ the two eigenvalues of the matrix above. We obtain
\begin{equation*}
  \lim_{\epsi\to 0,\infty} \frac{\alpha_1}{\alpha_G} = \lim_{\epsi\to 0,\infty} \frac{\alpha_2}{\alpha_G} = \infty
\end{equation*}
This enables us to conclude that, when the aspect ratio $\epsi$ tends either to zero
or infinity, the ratios above tends to infinity so that we cannot find a constant
$c>0$, independent of the aspect ratio $\epsi$, such that
\begin{equation*}
  c a_E(\vect{u}, \vect{u}) \leq s_E(\vect{u}, \vect{u}),
\end{equation*}
for all $\vect{u}\in\ker\projgrad$. It implies that the stabilization term is not
stable with respect to the aspect ratio.

\subsection{An alternative choice of the stabilization scaling}

Instead of using $\alpha_G$, let us use
\begin{equation}
  \label{eq:alphaNdef}
  \alpha_N = \frac19{\abs{E}}\trace(\hat{\mat{C}})\trace(\inv(\mat{N}_c^\trans\mat{N}_c)).
\end{equation}
Both $\alpha_N$ and $\alpha_G$ are invariant with respect to rotation. Because of the
coefficient $\frac19$, we have that if $\mat{N}_c^\trans\mat{N}_c$ \textit{were}
diagonal with constant coefficient, then $\alpha_N$ and $\alpha_G$ would be
equal. But in general they differ and we have
\begin{equation}
  \label{eq:expalphan}
  \alpha_N = \frac{2\lambda + 6\mu}4\left(\epsi + \epsi^{-1} + \frac{2}{\epsi + \epsi^{-1}}\right)
\end{equation}
It implies that
\begin{equation*}
  \lim_{\epsi\to 0} \frac{\alpha_1}{\alpha_N} = \lim_{\epsi\to \infty} \frac{\alpha_2}{\alpha_N} = \frac{16}{3}\frac{\lambda + 2\mu}{\lambda + 3\mu}
\end{equation*}
and
\begin{equation*}
  \lim_{\epsi\to\infty} \frac{\alpha_1}{\alpha_N} = \lim_{\epsi\to 0} \frac{\alpha_2}{\alpha_N} = \frac{16}{3}\frac{\mu}{\lambda + 3\mu}.
\end{equation*}
Therefore, for this choice of $\alpha$, there exist two constants $c_1,c_2>0$ which
are \textit{independent} of the aspect ration $\epsi$ and such that
\begin{equation*}
  c_1 a_E(\vect{u}, \vect{u}) \leq s_E(\vect{u}, \vect{u}) \leq c_2 a_E(\vect{u}, \vect{u})
\end{equation*}
for all $\vect{u}\in\ker\projgrad$. We can conclude that the stabilization provided
by $\alpha_N$ is stable with respect to the aspect ratio, at least for
quadrilaterals. Let us now try to explain the motivation back the introduction of
$\alpha_N$. We denote by $\lambda_i$ the singular values of $N_c$ and introduce the
following averages
\begin{equation*}
  \lambda_{\text{arithm}} = (\sum_{i=1}^{d}\lambda_i^2)^{\frac12}\quad\text{ and }\quad \lambda_{\text{harm}} = (\sum_{i=1}^{d}\lambda_i^{\frac12})^{2},
\end{equation*}
which, for simplicity, we refer to as \textit{arithmetic} and \textit{harmonic}
averages. Note that the matrix $\mat{N}_c$, which is given in (26), accounts for the
geometry and the unit of each coefficient is a unit length. We could therefore
interpret the values of $\lambda_{\text{arihm}}$ and $\lambda_{\text{harm}}$ as
characteristic lengths of the cell. Using these values, we can rewrite the scaling
coefficients as
\begin{equation*}
  \alpha_G = \frac{1}{\lambda_{\text{arithm}}^2}\abs{E}\trace(\hat{\mat{C}})\quad\text{ and }\quad
  \alpha_N = \frac{1}{\lambda_{\text{harm}}^2}\frac{\abs{E}}{9}\trace(\hat{\mat{C}}),
\end{equation*}
so that the difference between the two scalings is that they consider different type
of averages. Let us use eigenmodes to estimate the energy in each direction. For
simplicity, we consider the Laplace equation and the normalized energy of the mode
$\phi_i(x) = \cos(\frac{\pi}{2h_i}x_i)$ in the $i$-th direction is given by
\begin{equation*}
  \frac{\int_K\abs{\grad\phi}^2\,dx}{\int_K\abs{\phi}^2\,dx} = \frac{\pi^2}{(2h_i)^2},
\end{equation*}
from which we infer that a typical scale for the energy in the direction $x_i$
is given by $\frac{1}{h_i^2}$. If we consider a linear combination of such
unidirectional functions and neglect the interactions between them, then we
are naturally led to consider the sum
\begin{equation*}
  \sum_{i=1}^d\frac1{h_i^2}
\end{equation*}
as a typical scale for the energy. To obtain a typical length, we end up by taking
the harmonic average as defined above.

\section{Numerical test cases}

The great advantage of VEM methods is that they are valid for very general grids
including non-convex cells and more than one face between two cells,
\cite{beirao2013basic}. This property can be used to avoid curved faces on general
cells, simply by triangulating the surface. The VEM theory does not cover curved
surfaces and in the next examples we investigate the need for triangulation in 3D.

\subsection{A two-dimensional compaction case}

\begin{figure}[h]
  \begin{center}
    \includegraphics[width=0.8\textwidth]{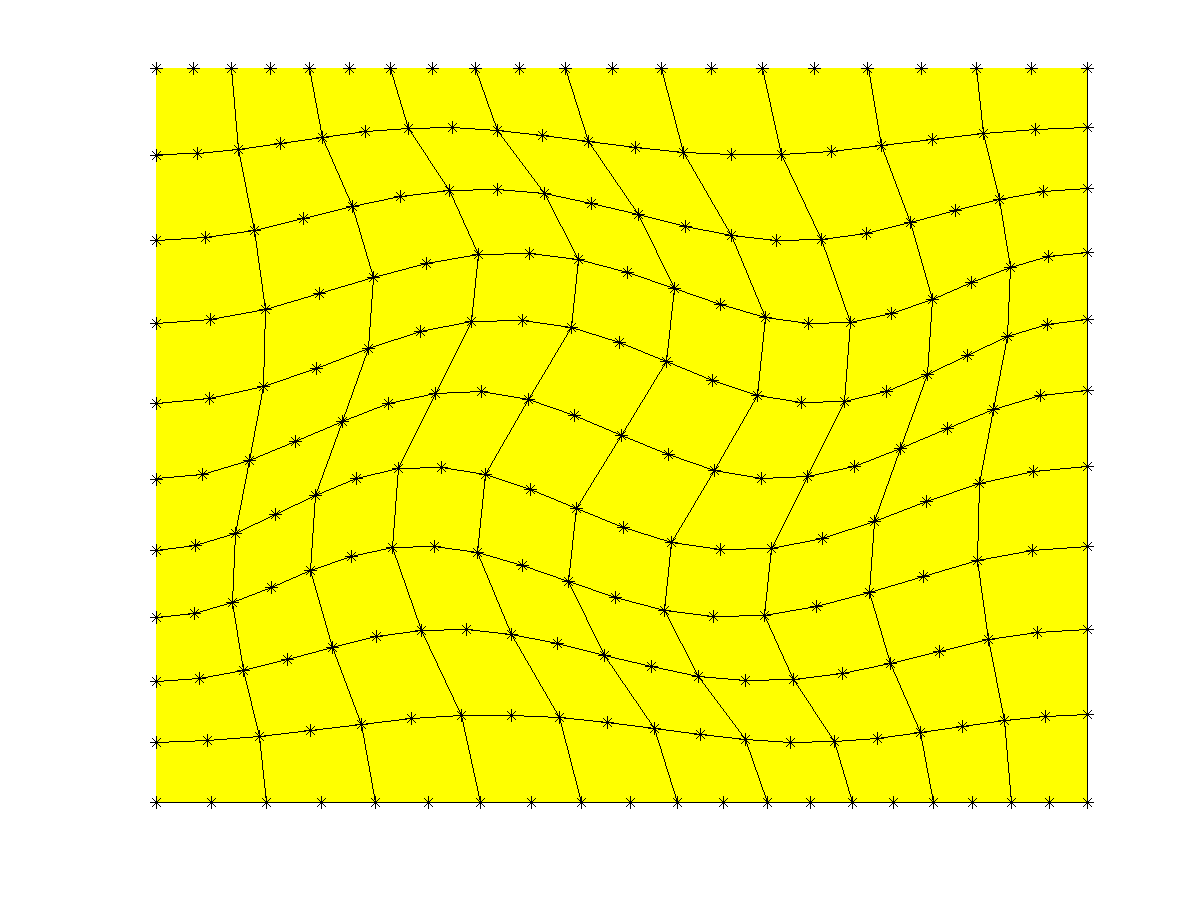}
    \caption{A twisted Cartesian grid is obtained by starting from a regular
      Cartesian grid and moving the nodes, here by using a smooth given displacement
      field. We plot the grid that is obtained after adding one
      extra node on each horizontal edge. Such grid is used to demonstrate the
      failure of the stabilization term where the aspect ratio is increased. The grid
      plotted here is the reference grid with aspect ratio, by definition, equal to
      $L_x/L_y=1$.}
  \end{center}
  \label{fig:2dcompactiongrid}
\end{figure}

\smallheading{Case description}: We consider a rectangular domain made of an
isotropic material with the following properties,
$\rho=\num{3e3}\si{\kilogram\per\cubic\meter}$, $E = \num{3e8}\si{\pascal}$ and
$\nu=0.3$. The vertical length of the grid is $L_y = \num{15}\si{\meter}$ and the
horizontal length will by determined by the aspect ratio $L_y/L_x$. Different values
of the aspect ratios will be tested. The boundary conditions are zero displacement at
the bottom, rolling boundary condition on the sides, that is no displacement in the
normal direction and no force in the tangential direction. At the top, we have no
force and free displacement. Even if the model is two-dimensional, we have to set up
boundary conditions for the third dimension, perpendicular to the plane, as the
material is going to expand or withdraw in this direction due to the Poisson
ratio. We impose zero displacement in the perpendicular direction, the other standard
choice being no force in that direction. The load term is gravitation, that is, a
constant vertical vector pointing downwards and we simulate the situation where the
material is going to subside by the effect of its own weight, hence the name of
\textit{compaction}. An analytical solution is available for this case and given by
\begin{equation}
  \vect{u}=[0, \gamma (L_y^2 - (y-L_y)^2)] \quad \gamma= \frac{g \rho}{2 C_{2,2}},
\end{equation}
where $C_{2,2}$ is the second diagonal coefficient of the stiffness matrix $\mat{C}$.
We start with a Cartesian grid that we twist in order to avoid artifact effects from
symmetries. We will refer to this grid as the \textit{twisted grid}. We consider a
variation of this grid where we add extra degrees of freedom in the form of extra
nodes on the horizontal edges, see Figure \ref{fig:2dcompactiongrid}. The motivation
for introducing such extra nodes is explained in the next paragraph.

\smallheading{Results}: We test the three different implementations of the load term,
as described in Section \ref{subsec:loadterm}. It is important to note that the first
alternative, which uses the projection operator, see the definition in
\eqref{eq:exactf1}, is \textit{exact} in the case we are considering. Indeed, since
$\vect{f}=\rho\vect{g}$ is constant, we have
\begin{equation*}
  \proj_i^0(\vect{f}) =
  \vect{f} 
\end{equation*}
so that no error is introduced by the assembly of the load term. In the remaining, we
will refer to this implementation of the load as the \textit{exact} load term. In
comparison, the third method is not exact, as the potential function, here given by
$\psi = \rho g y$, is approximated by a cell-valued function. For the stabilization
term, we test the two scaling variables $\alpha_G$ and $\alpha_N$ presented in
Section \ref{subsec:stabaspectrat}.

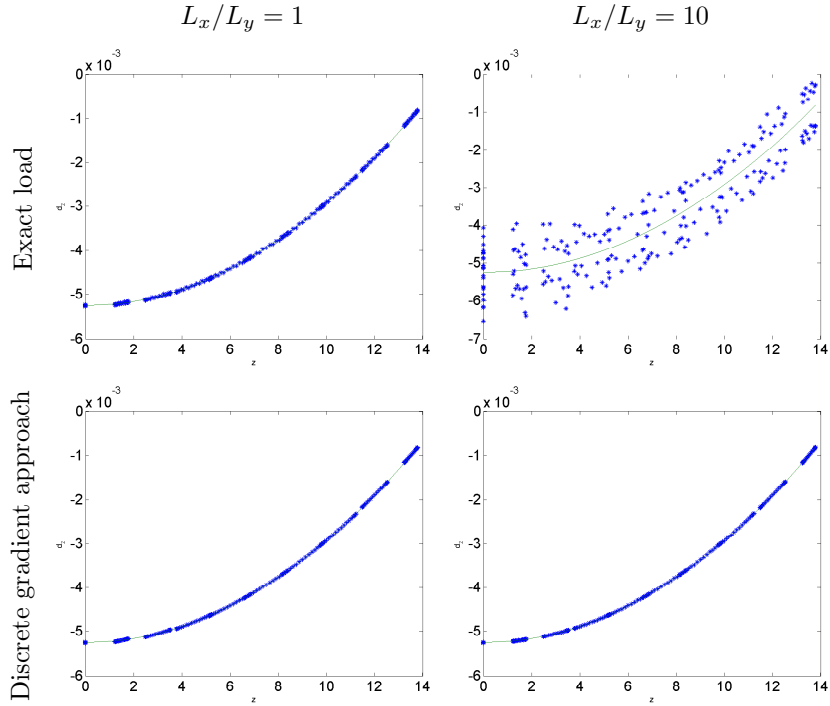
\begin{figure}[h]
  \begin{center}
    \begin{tikzpicture}
      \matrix at (0,0) {    
        &\node  {$L_x/L_y=1$}; &  \node {$L_x/L_y=10$};\\
        \node[transform shape, rotate = 90] {Exact load}; &\node{\includegraphics[width=0.4\textwidth]{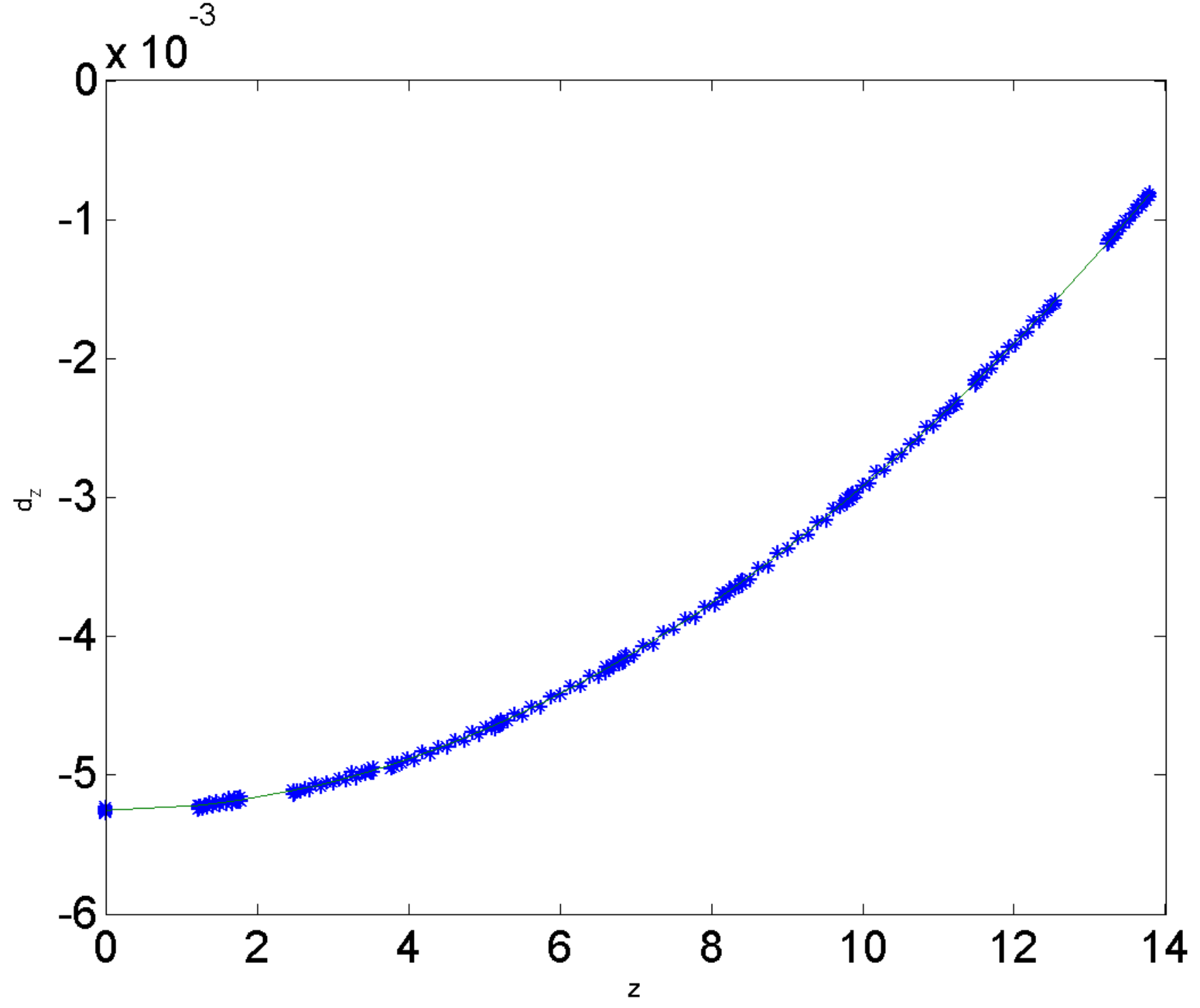}};
        &\node{\includegraphics[width=0.4\textwidth]{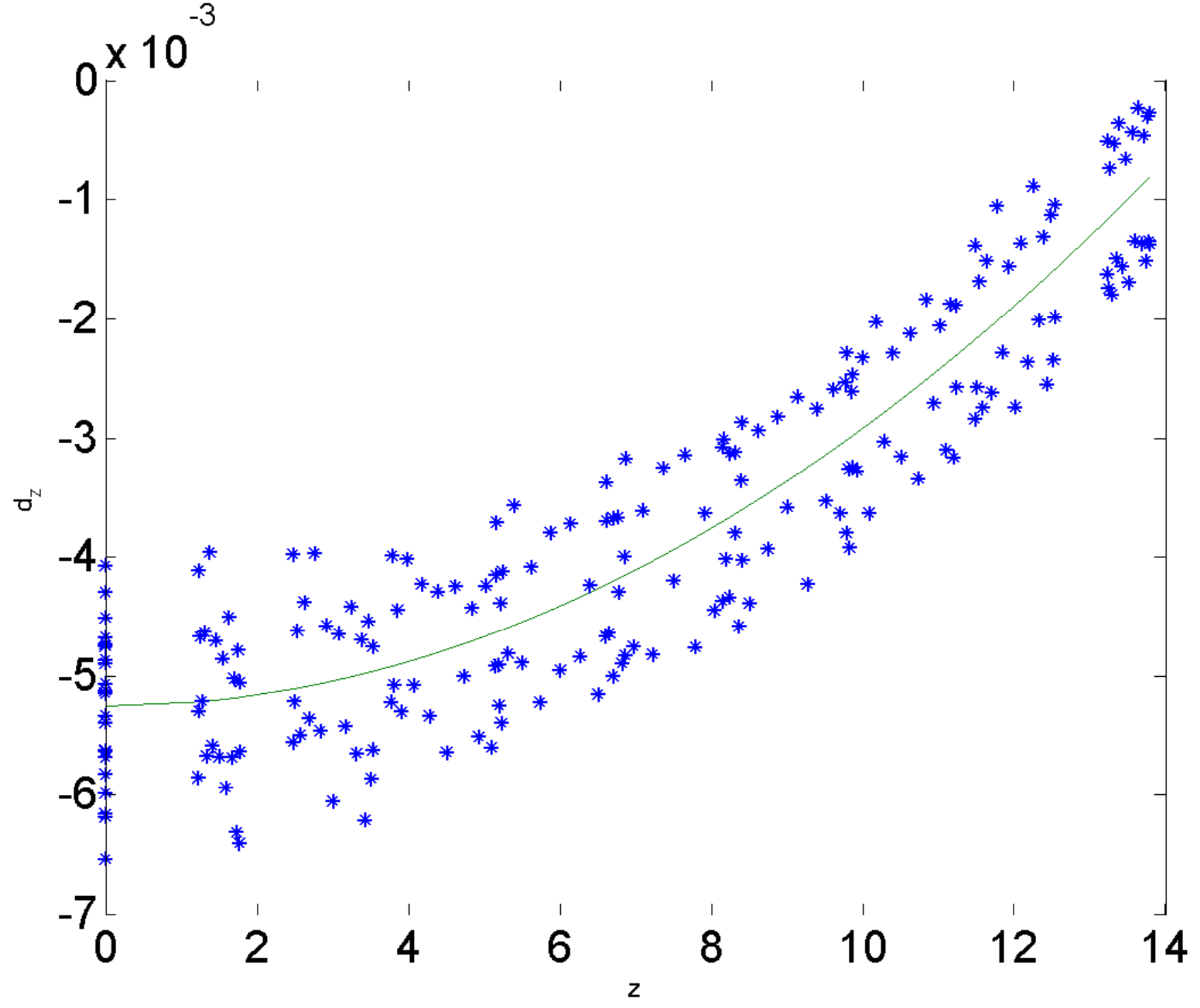}};\\
        \node[transform shape, rotate = 90] {Discrete gradient approach};&\node{\includegraphics[width=0.4\textwidth]{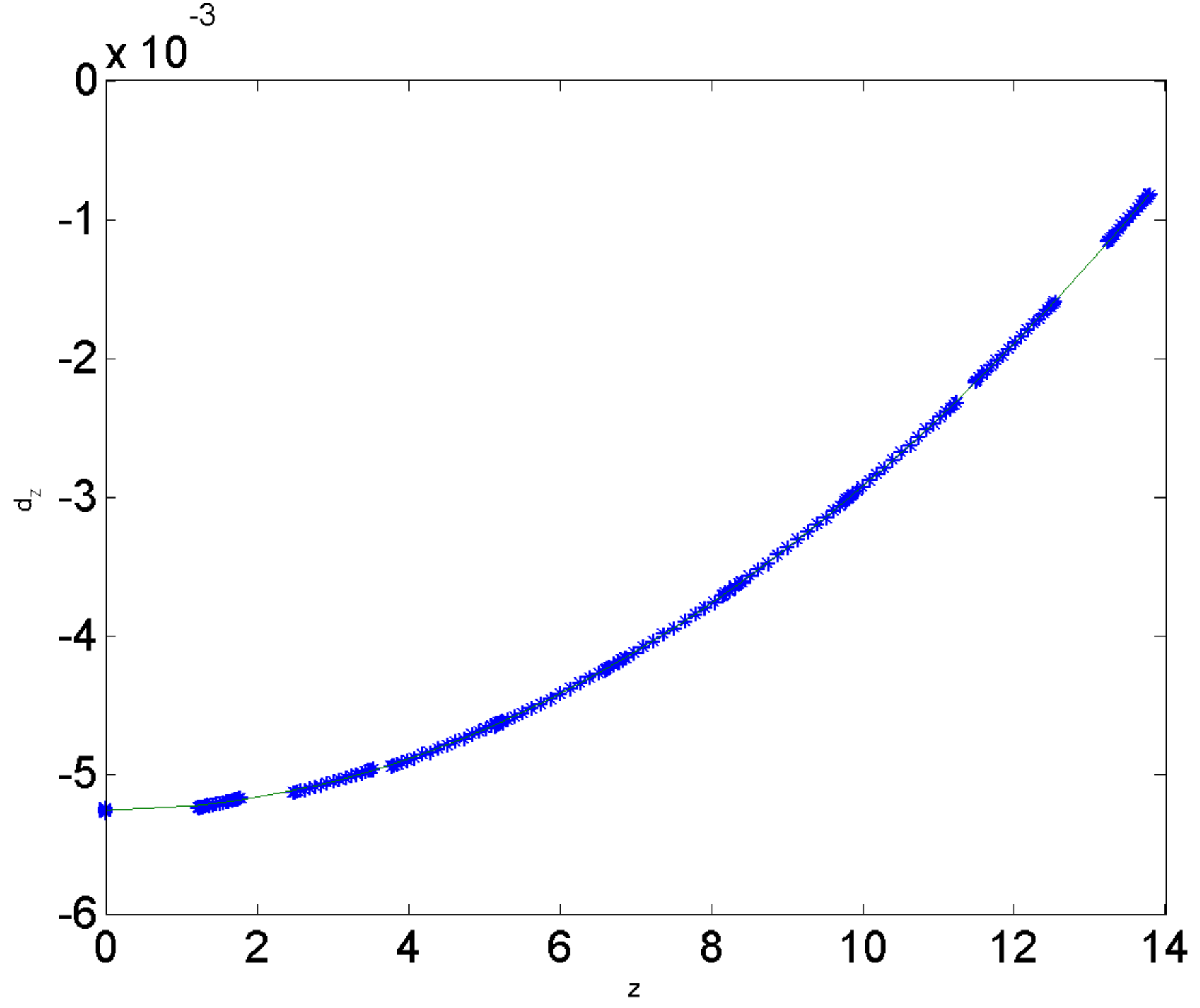}};
        &\node{\includegraphics[width=0.4\textwidth]{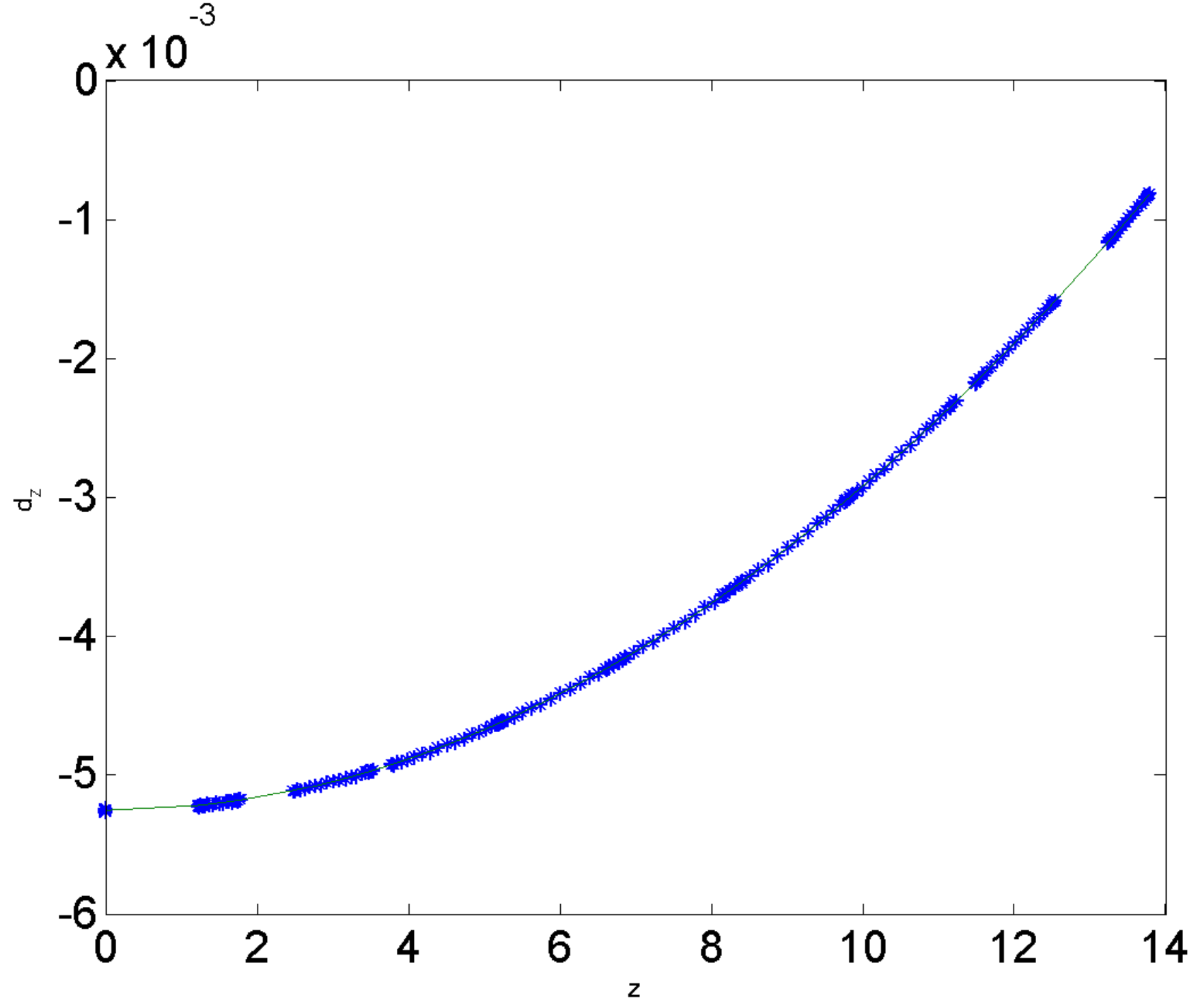}};\\
      };    
    \end{tikzpicture}
  \end{center}
  \caption{We plot the computed displacement in the vertical direction for the 2D
    compaction example. The result for every node of the grid is represented as a
    dot, where the $x$-coordinate of the node corresponds to the vertical position of
    the node and the $y$-coordinate corresponds to the value of the vertical
    displacement computed at the node. The analytical solution is plotted as a
    continuous line. For these plots, the twisted Cartesian grid with an extra node
    on each horizontal edges, see Figure \ref{fig:2dcompactiongrid}, has been
    used. The left column is for aspect ratio $1$ and the right is for aspect ratio
    $10$. For the first row, the exact load calculation based on the exact
    integration of the VEM basis function has been used, while the lower row
    corresponds to the discrete gradient approach. We use the scaling factor
    $\alpha_G$ as proposed in \cite{gain2014}, see \eqref{eq:defalphagain}. We
    observe that, for the exact load computation, the solution quickly deteriorates
    when the aspect ratio is increased while the results remain good for the discrete
    gradient approach.}
  \label{fig:2D_compaction}
\end{figure}

We start with the scaling variable $\alpha_G$ taken from \cite{gain2014} and the
exact load implementation. We use the grid with extra nodes. For such grid, each cell
gets two extra degrees of freedom. However, these extra degrees of freedom do not
enrich the approximation space as they do in the case of a finite element method. The
VEM method retains the same degree of accuracy, that is first order in our case. The
extra basis functions introduced by the extra degrees of freedom are handled by the
stabilization term. But the stabilization term only guarantees that these extra
functions do not break the ellipticity of the system but it is an artificial term
which cannot add any accuracy. Therefore, by adding an extra node on the edges, we
increase the relative importance of the stabilization term, so that its deficiency
will be more apparent. As predicted by the results of Section
\ref{subsec:stabaspectrat}, we observe a severe dependence on the aspect ratio. When
the aspect ratio is minimal, that is $L_x/L_y=1$, then the solution is close to the
analytical one but, when the aspect ratio is increased to $L_x/L_y=10$, by stretching
the grid in the horizontal direction, the results deteriorate severely, see the top
panels in Figure \ref{fig:2D_compaction}. We run the same simulations but, instead of
the exact load term, we use the load term computed by the discrete gradient
operator. Then, the results do not deteriorate as the aspect ratio is increased.

\begin{figure}
  \begin{tabular}{c c}
    \includegraphics[width=0.45\textwidth]{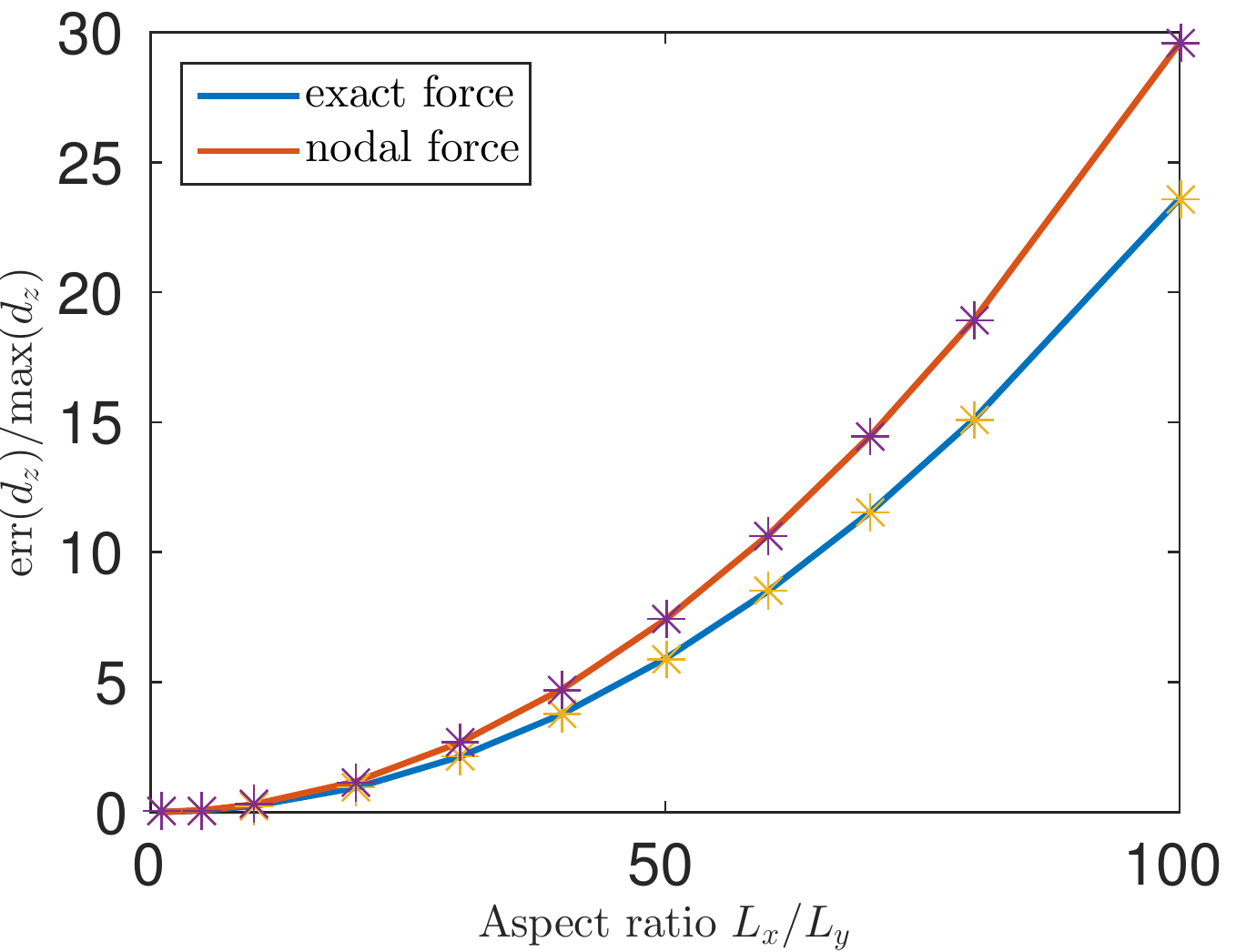}&
    \includegraphics[width=0.5\textwidth]{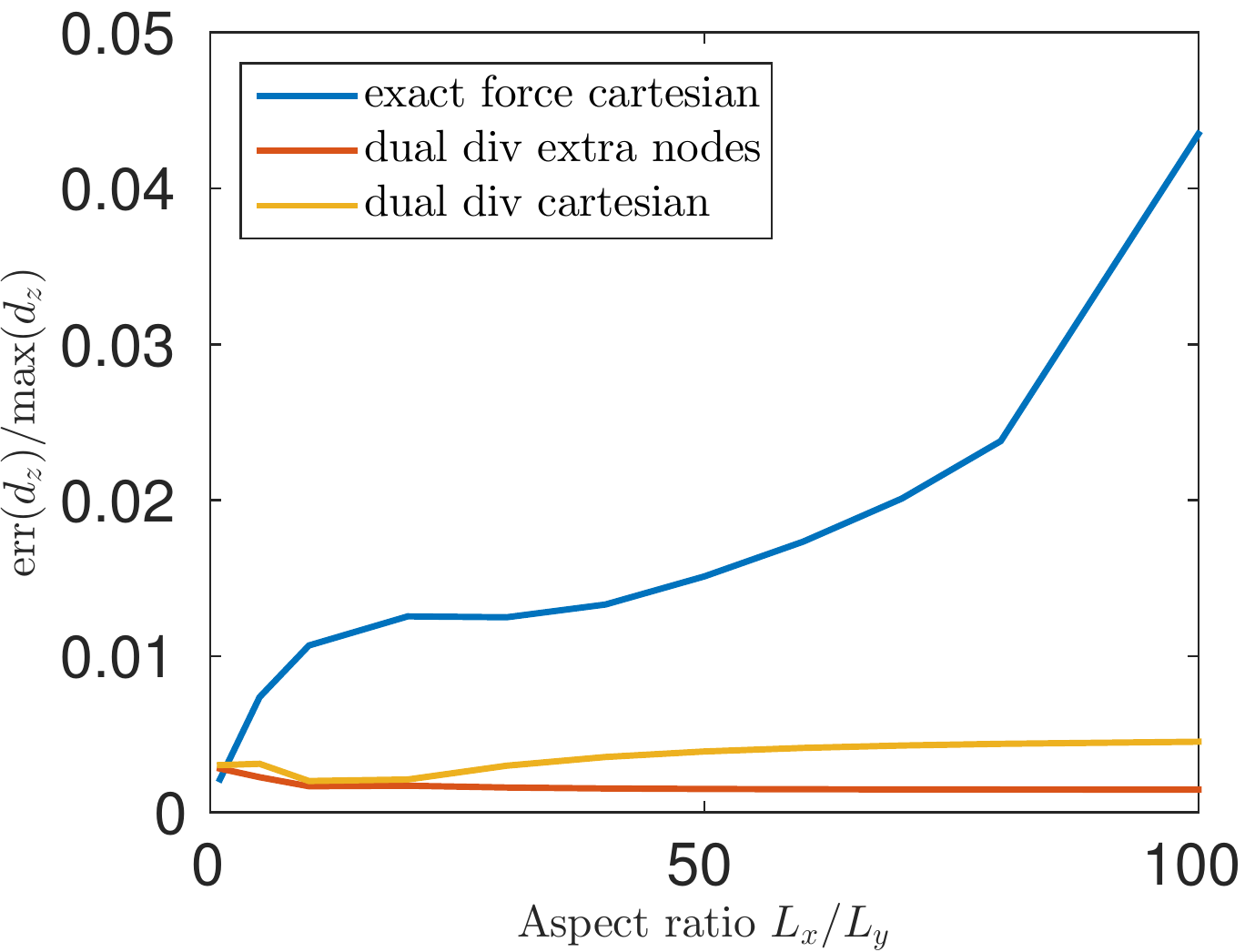}
  \end{tabular}
  \caption{Plots of the maximum error in the vertical displacement as a function of
    the grid aspect ratio. The left figure shows results for the exact integration
    method in the case with extra nodes on horizontal faces. We observe that the
    method fails as the error blows up. The extra points in this plot are reference
    points that indicate a quadratic scaling of the error with respect aspect
    ratio. In the right figure, the results are shown for the exact method on the
    twisted grid without extra nodes and for the discrete gradient approach on the
    same grid with and without the extra nodes.}
  \label{fig:2D_compaction_as}
\end{figure}

In Figure \ref{fig:2D_compaction_as}, we plot the error in displacement as a function
of the aspect ratio (from 1 to 100) for the different grid cases and the three
implementation of the load term.  The left figure shows that the exact load and nodal
load calculations fail for the grid with extra nodes. The error apparently follows a
second order growth, that is $err\sim (L_y/L_x)^2$. The plot on the right shows the
error for the twisted grid without extra nodes for the exact load computation and the
error for both grids for the discrete gradient approach. All the methods give
reasonable results, but the exact load calculation seems to deteriorate more than the
others. The discrete gradient approach is stable in both cases. Note that, if we had
used a grid without disturbance, all the methods would give exact results for the
grid without extra nodes on the faces while the extra node case will still fail for
the exact load calculation. The reason is that, in the non disturbed case with no
extra nodes, all the implementations of the load term give the same result in the
case of a constant vertical load term.

\begin{figure}[h]
  \centering
  \begin{center}
    \includegraphics[width=0.6\textwidth]{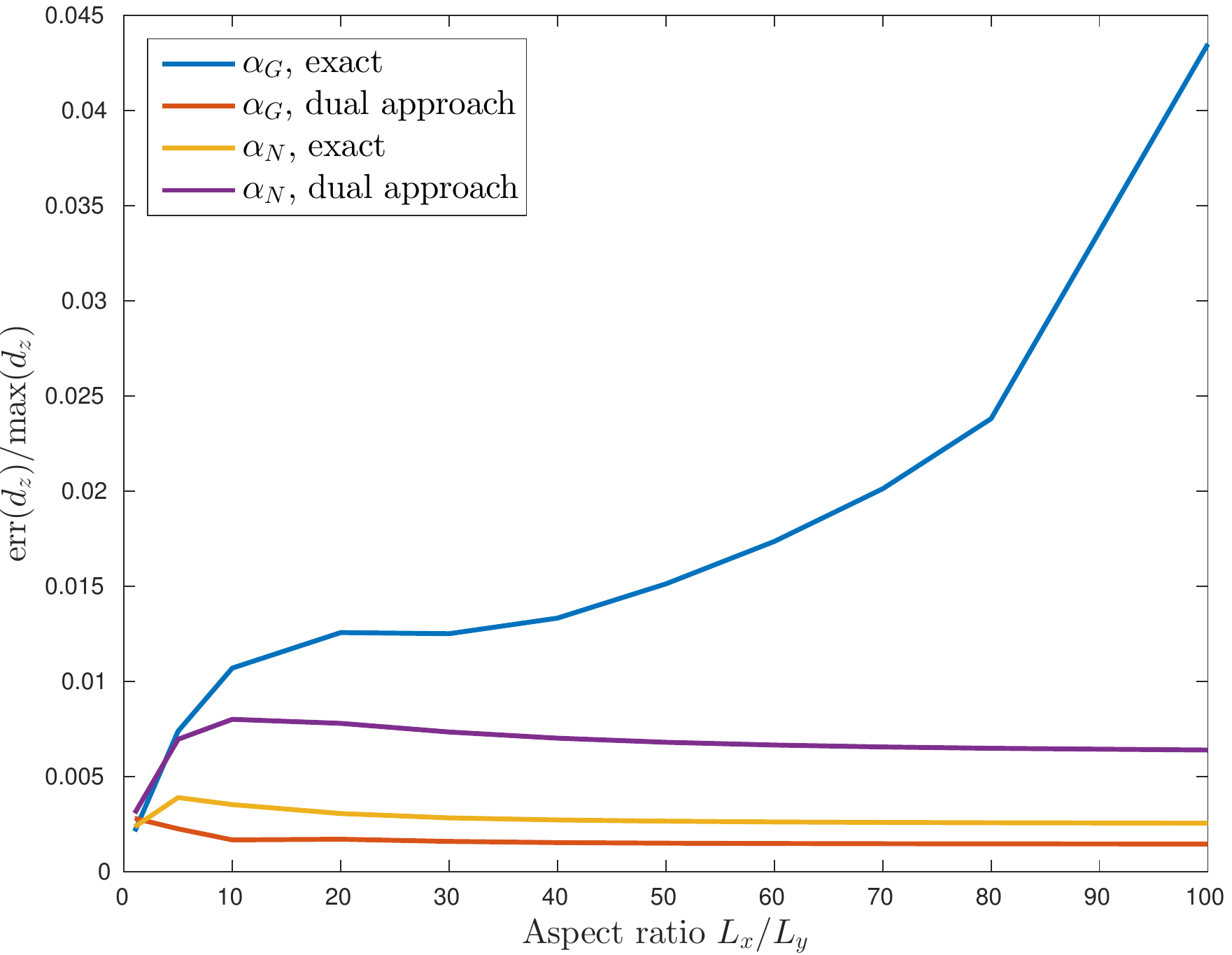}
  \end{center}
  \caption{Comparison between the two scaling constants $\alpha_G$ and $\alpha_N$. We
    plot the error of the vertical displacement, as the aspect ratios is
    increased. Here, we use the grid without the extra nodes. We observe that the
    scaling constant $\alpha_N$ yields stability with respect to aspect ratio,
    independently of which method is used to compute the load term.}
  \label{fig:err_as_dz_force_scaling}
\end{figure}

Finally in Figure \ref{fig:err_as_dz_force_scaling}, we consider the scaling
$\alpha_N$ introduced in \eqref{eq:alphaNdef}, which is stable with respect to aspect
ratio. The error does not grow as the aspect ratio is increased, as opposed to
$\alpha_G$. The use of $\alpha_N$ deteriorates the solution computed using the
discrete gradient approach, while it significantly improves the solution using the
exact method. However, this conclusion is difficult to extend to more general
cases. The value of $\alpha_N$ has been derived from an analysis done on regular
quadrilaterals and we observe that the stability properties extend to a twisted
Cartesian grid. However, separate studies would have to be done for more complicated
shapes and also in 3D, where the situation is expected to be more
complicated. Indeed, while in 2D the aspect ratio is described by a scalar quantity
namely $\epsi=\frac{\Delta x}{\Delta y}$, in 3D we need 2 values, say $\frac{\Delta
  x}{\Delta z}$ and $\frac{\Delta y}{\Delta z}$, the third quantity $\frac{\Delta
  x}{\Delta y}$ being imposed by the fact that we will anyway require isotropic
stability. It means that a scalar approximation of the stabilization term, as given
in \eqref{eq:stabgain} and also in \cite{da2013virtual}, will not be enough. This
problem was noticed in \cite{Andersen15:ecmorxv}, and the exact stabilization term
corresponding to finite element was used there to study a poro-elastic response
function in the 3D case.

\smallheading{Comment}: We do not really understand why the discrete gradient
approach (Alternative 2) performs significantly better than the projection approach
which is exact in this case (Alternative 1). However, we note some fundamental
differences between the force-based methods (Alternatives 1 and 2) and the discrete
gradient approach, which may help to understand the differences in the results. As
explained in the previous section, see \eqref{eq:nodalforce}, the difference between
the methods is in how they divide the weights between the nodes. All the force-based
methods divide forces according to a weight for each node associated with volume
integrals. These weights are equal for \textit{all} Cartesian directions. In
contrast, the discrete gradient method uses weights associated with surface
integrals, so that the weights can depend on the direction, and the corresponding
degrees of freedom. These weights can be associated with the projected area of the
faces associated with a node divided by the projection of the cell in the same
direction. This is most easily seen from the expression in equation
\eqref{eq:local_faceforce}. In the case of the extra nodes on the edges, these nodes
will have associated weights in the horizontal direction only due to the tilt of the
grid and the weights will in the simple case be doubled of the corner nodes while the
exact case will give all nodes the same weights.  In \cite{gain2014}, the method
using node quadrature (Alternative 2) is considered, this will in the above case give
a smaller weight to the midpoint and behave worse for the case with extra node, as
seen in the left panel of Figure \ref{fig:2D_compaction_as}.

\clearpage

\subsection{Compaction 3D}
\label{subsec:compaction3D}
\begin{figure}
  \begin{tikzpicture}
    \matrix at (0,0) {
       \node {\includegraphics[width=0.3\textwidth]{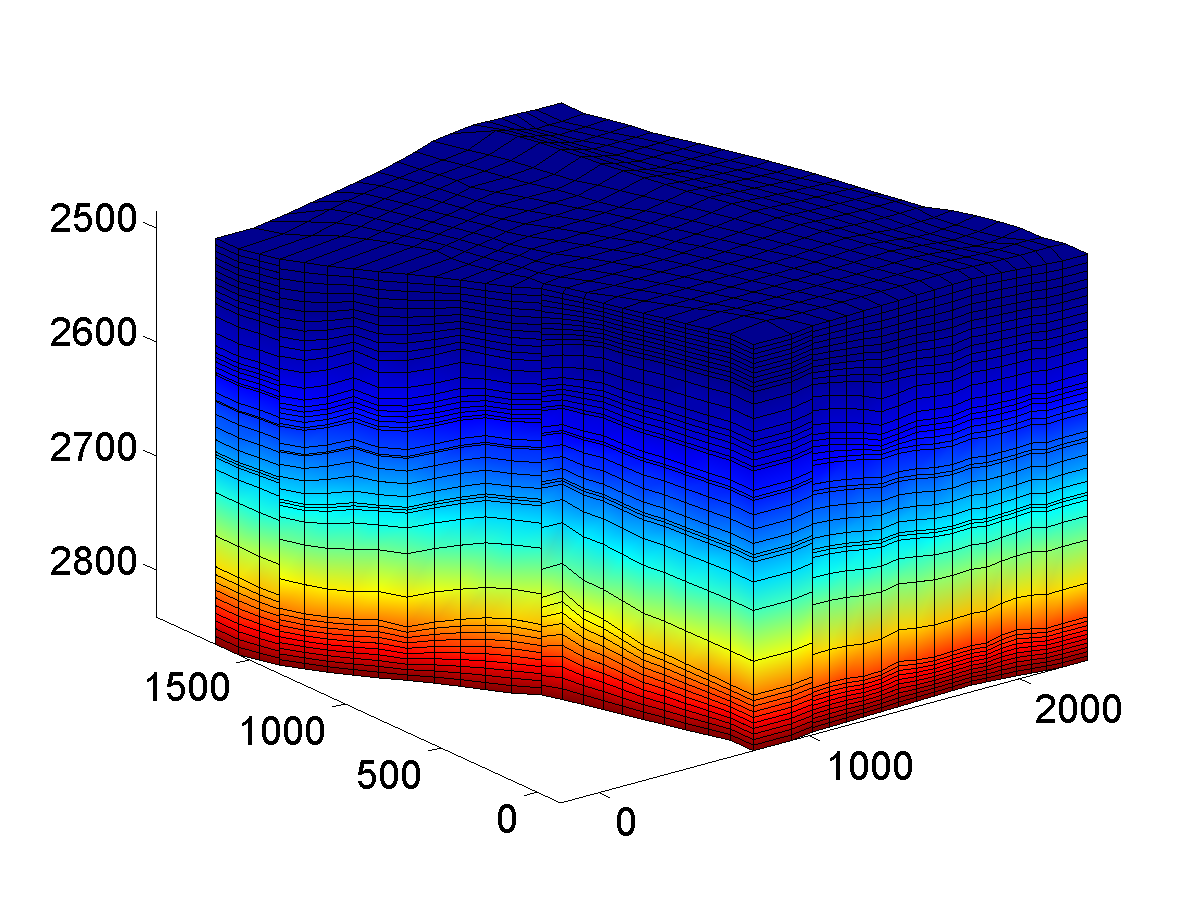}};
      &\node  (orig) {\includegraphics[width=0.3\textwidth]{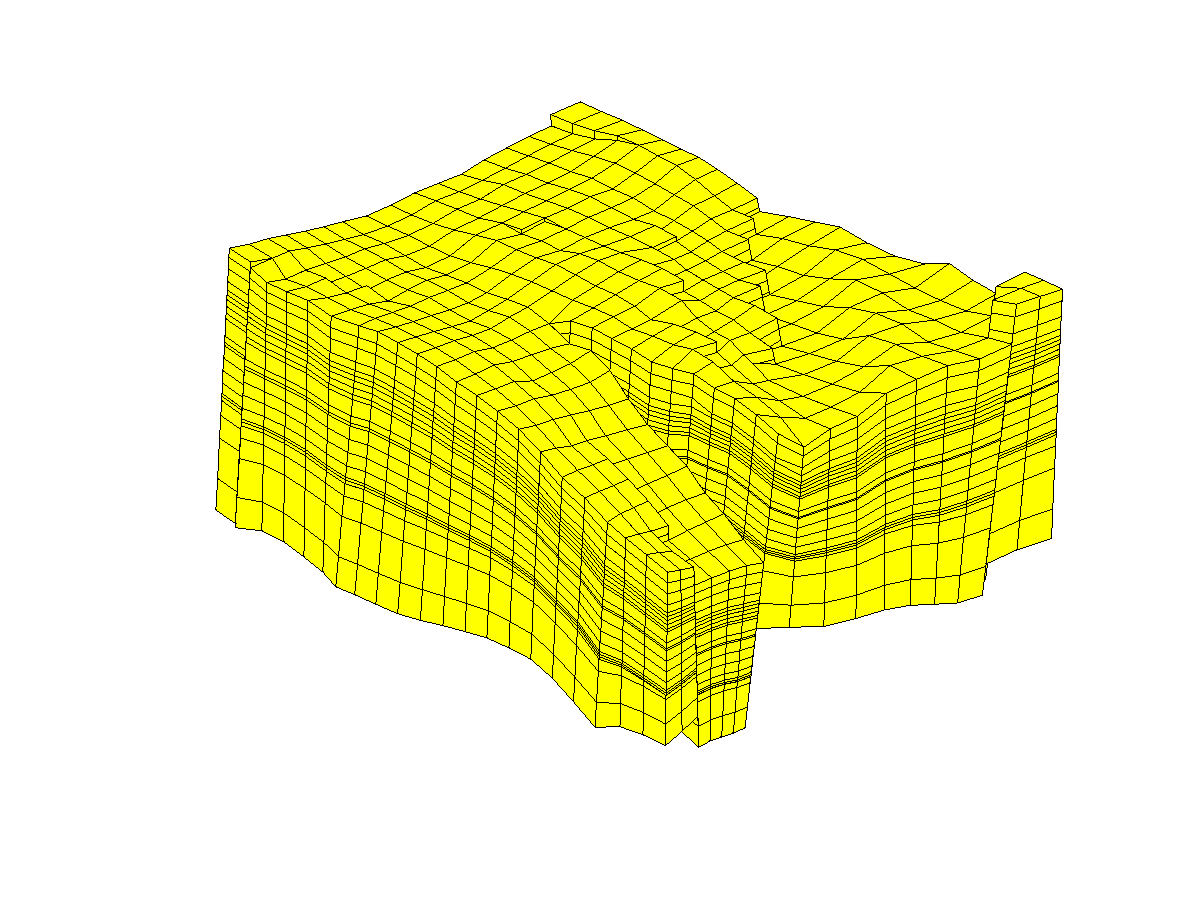}}; &
       \node  (vert) {\includegraphics[width=0.3\textwidth]{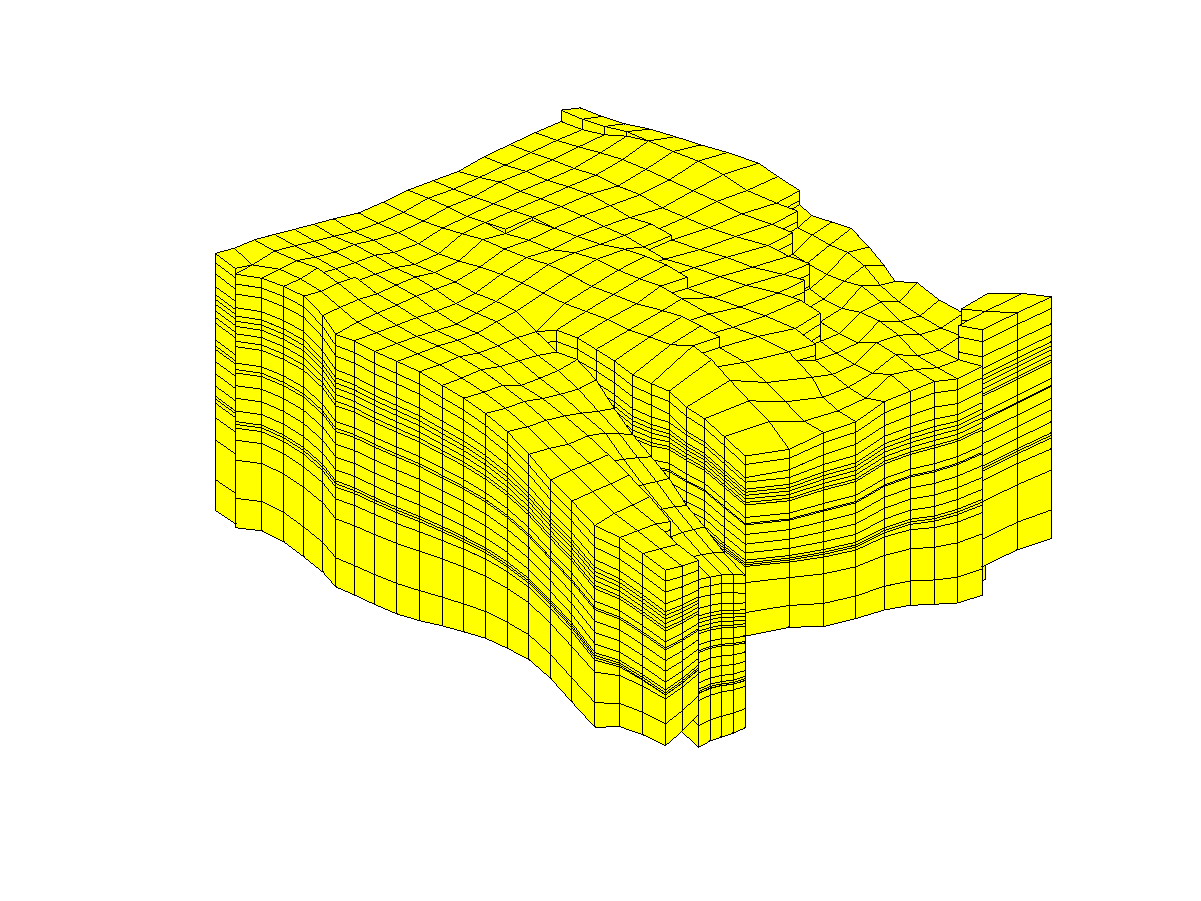}}; \\
      \node{\includegraphics[width=0.3\textwidth]{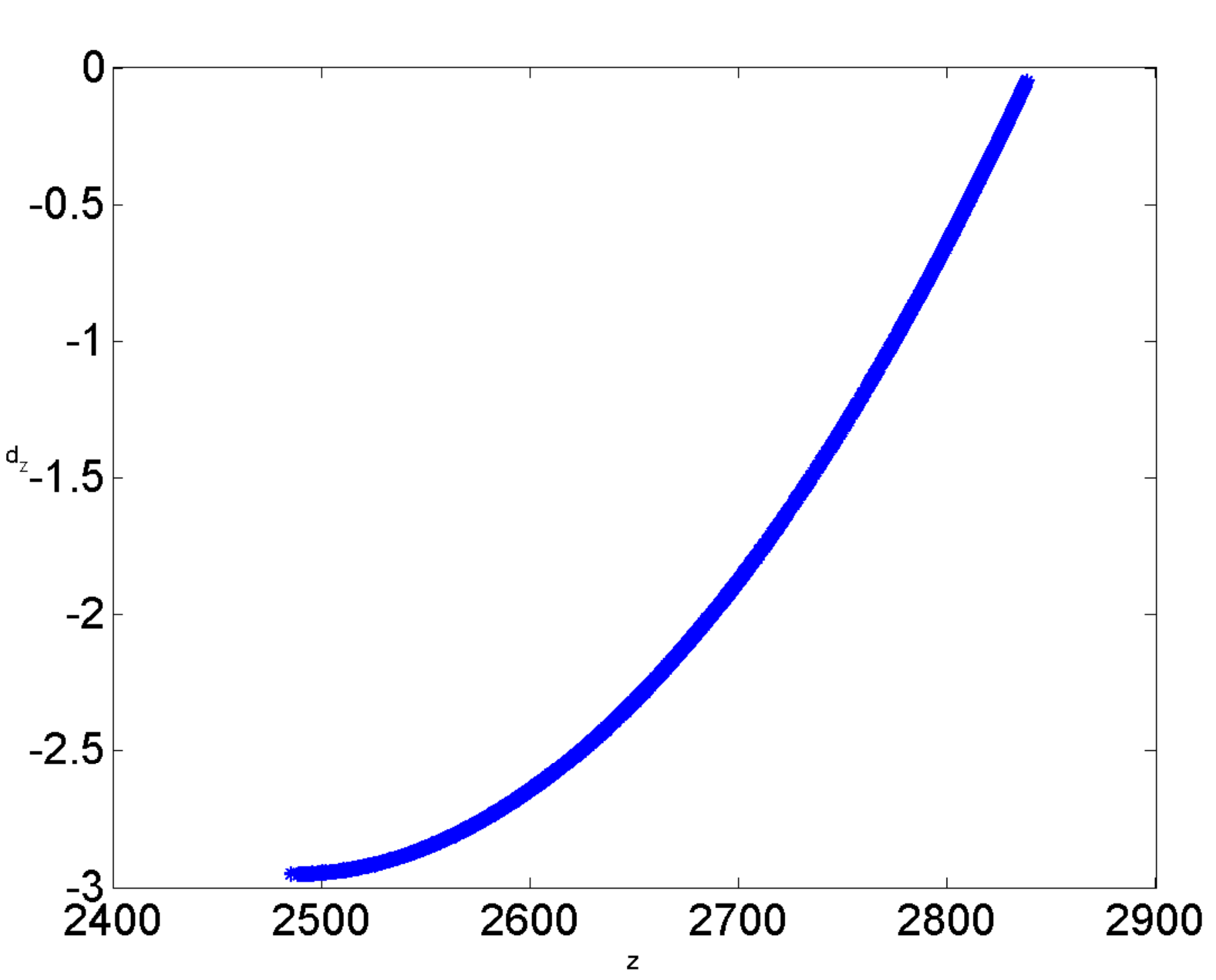}};&
      \node {\includegraphics[width=0.3\textwidth]{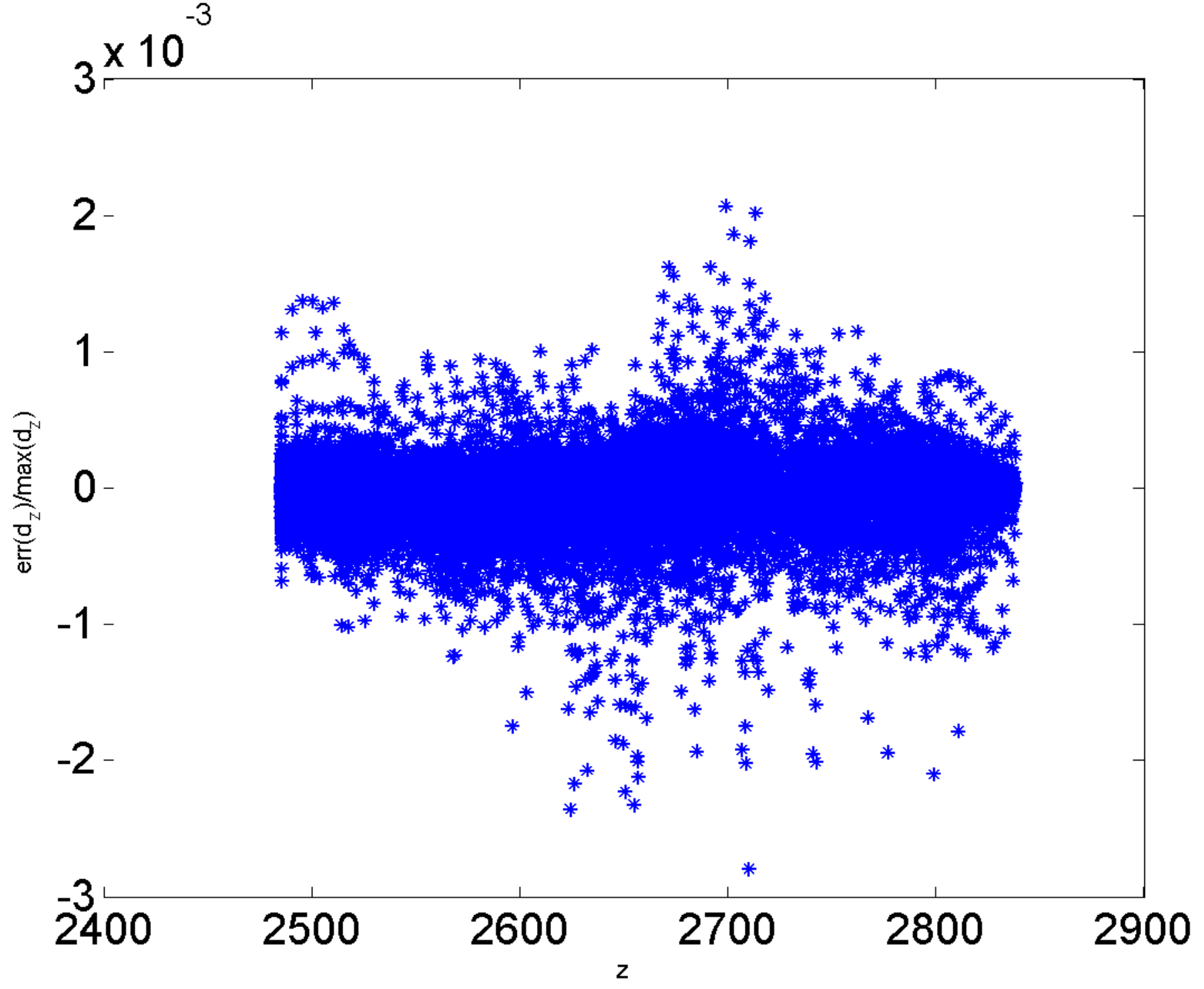}};&
      \node{\includegraphics[width=0.3\textwidth]{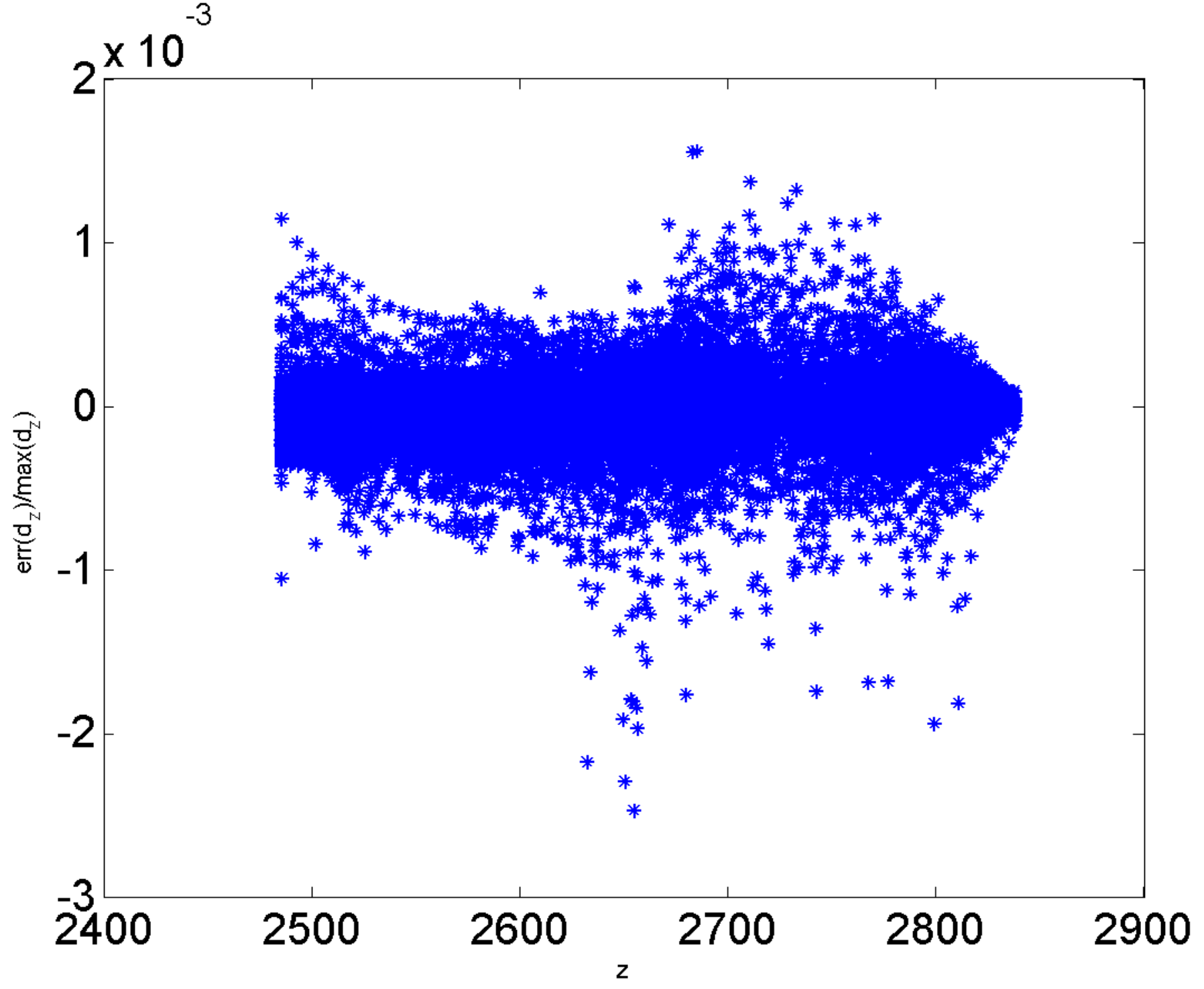}};\\
    };
    \node[anchor = south] at (orig.north) {Original pillars};
    \node[anchor = south] at (vert.north) {Vertical pillars};
  \end{tikzpicture}
  \caption{Effect of compression on a part of the Norne model. The first column shows
    the plot the vertical deformation on the grid (left), the original grid after
    removing the padding (middle) and the straightened grid where the pillars are
    made vertical (right). The figure in the upper left corner shows the bounding box
    which is used for the calculation, while the two other grids show the embedded
    Norne grid. The second row show the results for pure gravitational
    compression. The first column shows the vertical displacement while the second
    and third show the errors in the vertical displacement for the original and
    triangulated grid.}
  \label{fig:norne}
\end{figure}
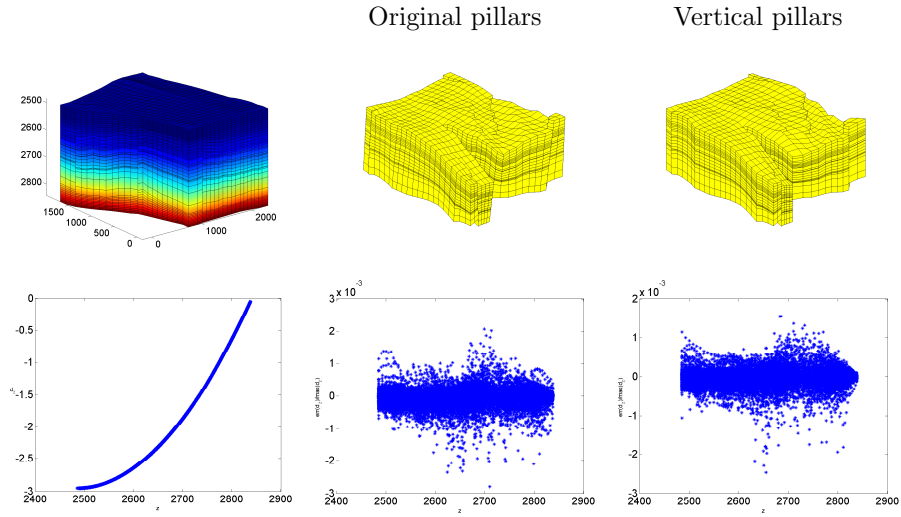

In order to investigate the performance of the VEM method on real reservoir
geometries, we use two grids which includes standard features of subsurface
models. The first one is based on a local sedimentary model called \sbed. The model
was used for upscaling permeability. Our version is $15m\times15\times3m$ with
logical Cartesian dimensions $15\times15\times333$.  The grid reflects two of the
basic properties of a sedimentary process, which are the layering and erosion
processes. For this type of grid, the challenge is the degenerate cells and the large
aspect ratios. The second model that we consider is taken from the open reservoir
model of Norne. The data for this model is freely available in the open dataset of
the Open Porous Media initiative \cite{OPM:webdata}. We extract a part of this model,
pad it on all sides to embed it in a regular prism, so that we can simply impose side
boundary conditions and directly compare the solution with the analytic solution of a
pure gravitational compression. The final full model and the embedded model with
faults are shown in Figure \ref{fig:norne}.

\begin{figure}[h]
  \centering
  \begin{tabular}[h]{c@{\hspace*{1cm}}c}
    \includegraphics[width=0.3\textwidth]{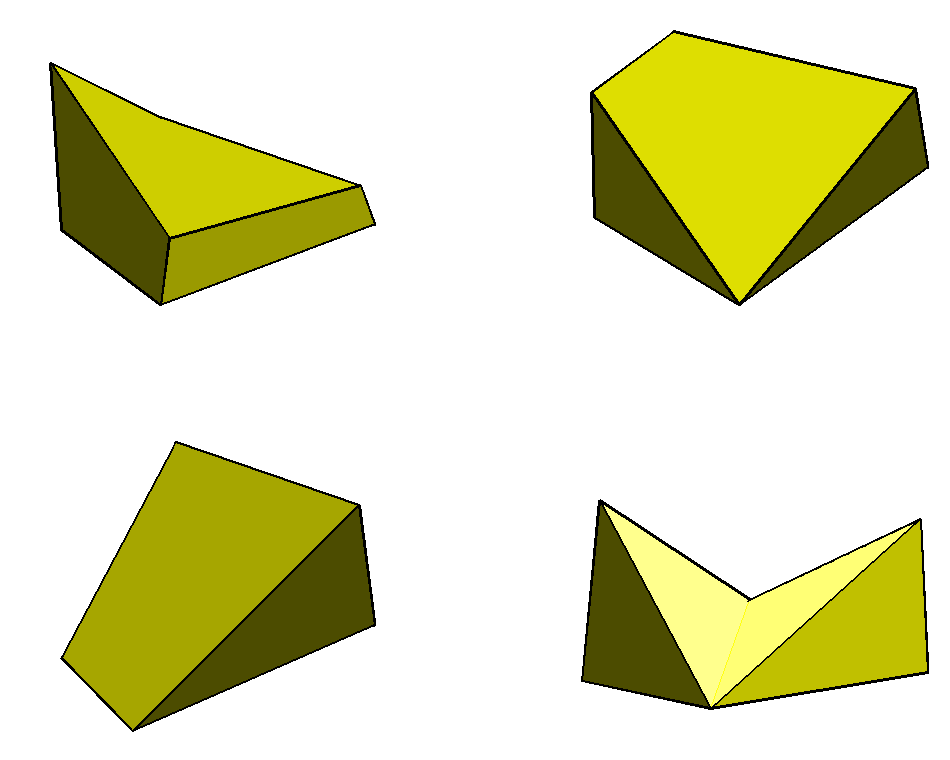}&\includegraphics[width=0.3\textwidth]{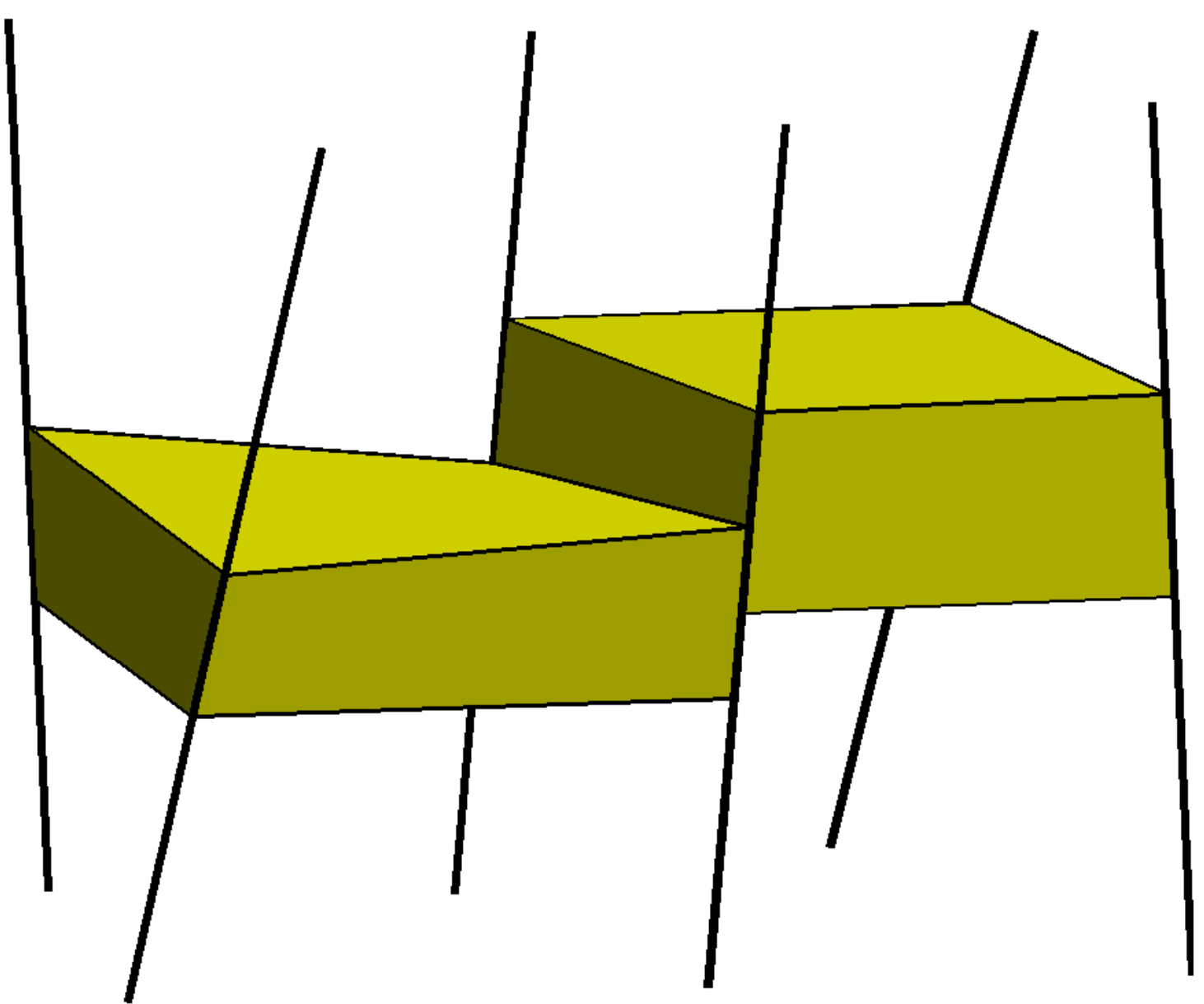}
  \end{tabular}
  \caption{On the right, two neighboring cells of a corner-point grid. On the left,
    examples of the irregular cell shapes that the corner-point format can produce.}
  \label{fig:cornerpointcells}
\end{figure}

Both models use a corner-point grids, which is a standard in the industry.  A
corner-point grid has an underlying two dimensional structure which is used to index
the pillars. Let us denote by $p_{i,j}$ and $q_{i,j}$ the bottom and top and the
pillar that is indexed by $(i,j)$. For each region contained between the four pillars
$(i,j)$, $(i+1.j)$, $(i+1,j+1)$ and $(i,j)$, points are defined on each of this
pillar in equal number. We denote those points by $x_{i',j'}^k$ for $i'\in\{i,i+1\}$
and $j'\in\{j,j+1\}$. Then, the region between the four pillars is meshed with
hexahedrons with eight corner points given such as $x_{i',j'}^{k'}$ for
$i'\in\{i,i+1\}$, $j'\in\{j,j+1\}$ and $k'\in\{k,k+1\}$.  This construction naturally
leads to irregular cell shapes and faces that are not planar, see the illustrations
given in Figure \ref{fig:cornerpointcells} .
Therefore we end up outside the theoretical framework of the VEM method, which only
cover planar polygonal faces. However, the computation of the stiffness matrix for
VEM relies on geometrical properties that are all available, either as exact or
approximated values (such as face areas, face normals, etc.), in the case of a
corner-point grid, so that the stiffness matrix can be assembled and a solution
computed. To evaluate the error that is introduced by this geometrical approximation,
we compare the solution obtained this way with the solution that is obtained after
triangulating the non planar surfaces, by adding a point in the middle of the
faces. For such grid, the faces will be planar and the theoretical framework of the
VEM method applies.

\begin{figure}
  \begin{tikzpicture}
    \matrix at (0, 0) {
      \node{\includegraphics[width=0.3\textwidth]{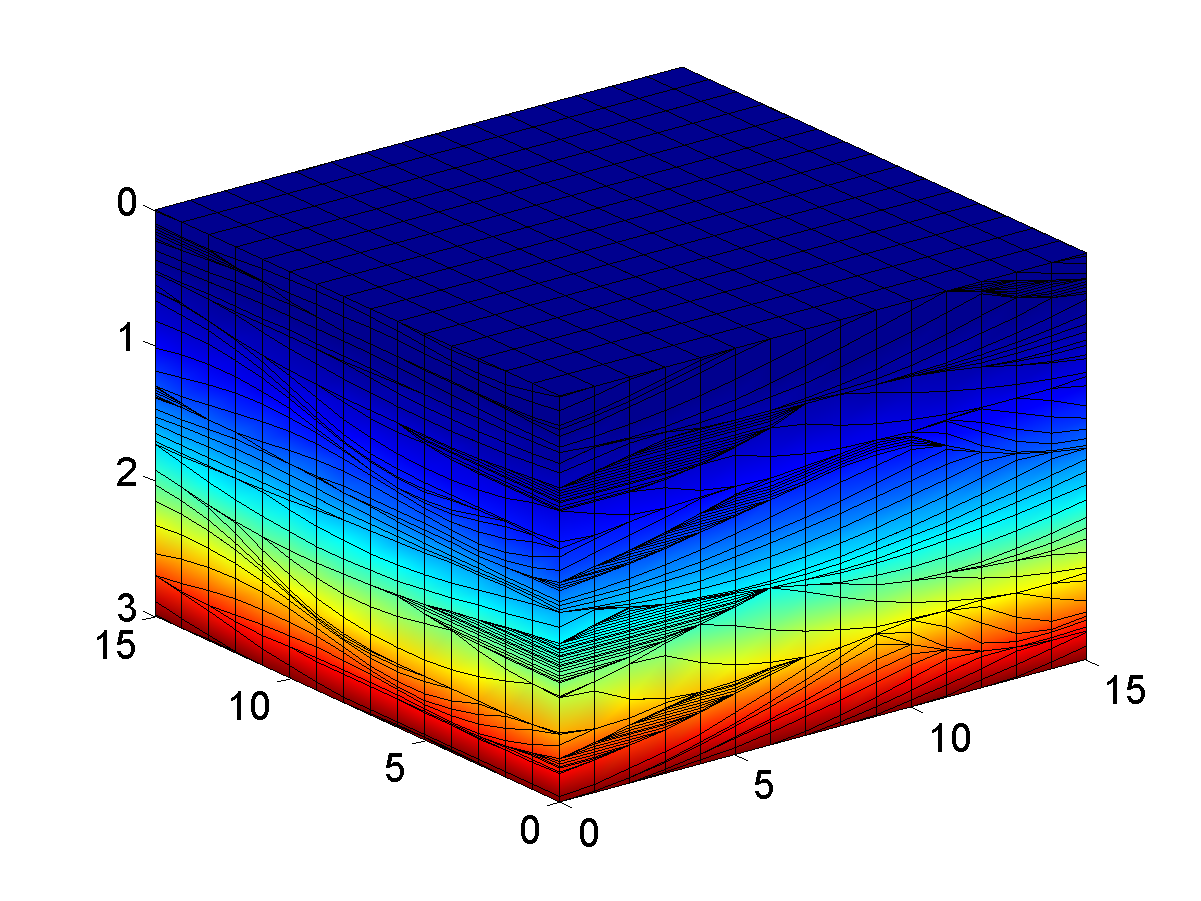}};&
      \node (notri) {\includegraphics[width=0.3\textwidth]{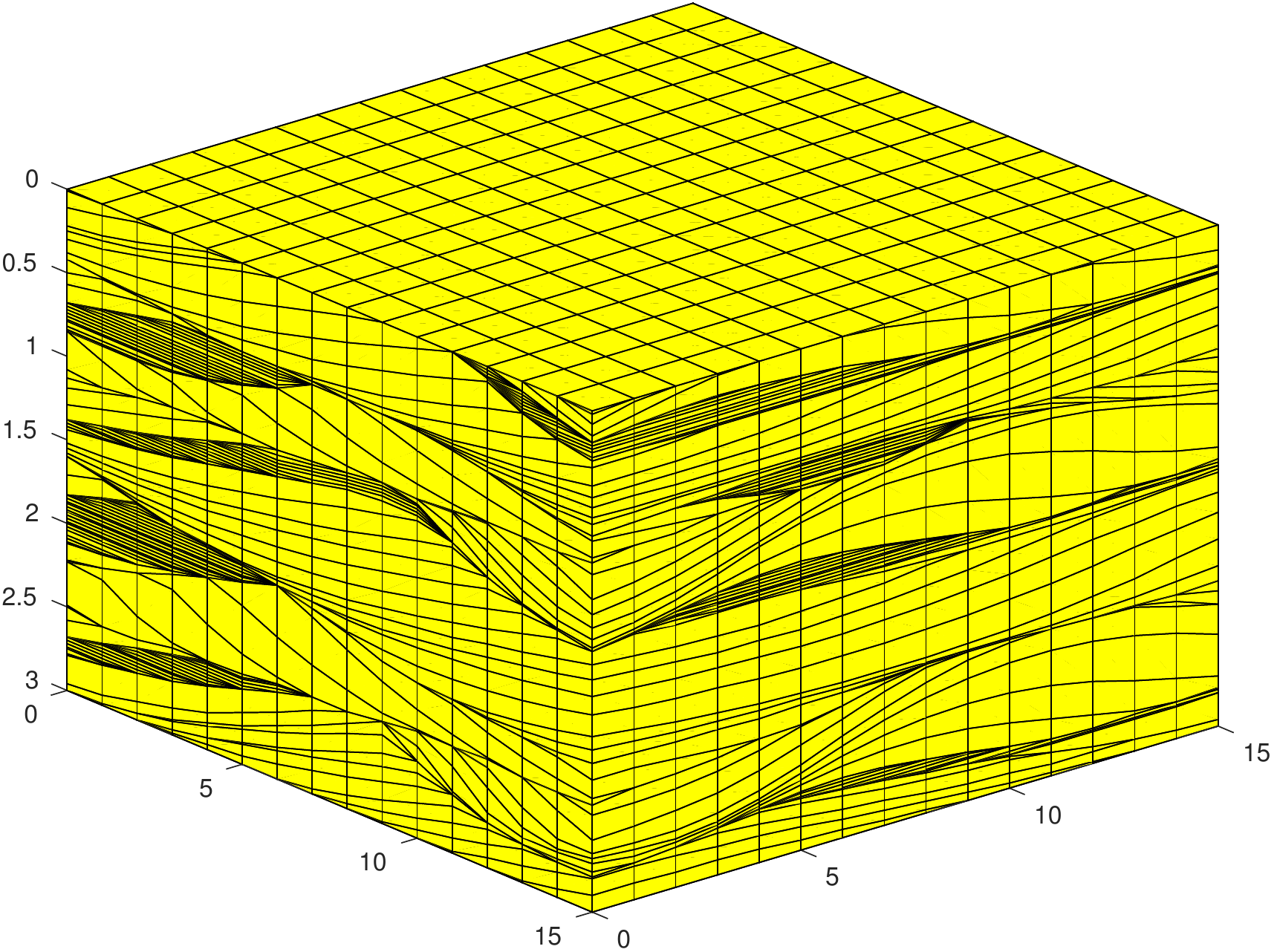}};
      &            
      \node (tri) {\includegraphics[width=0.3\textwidth]{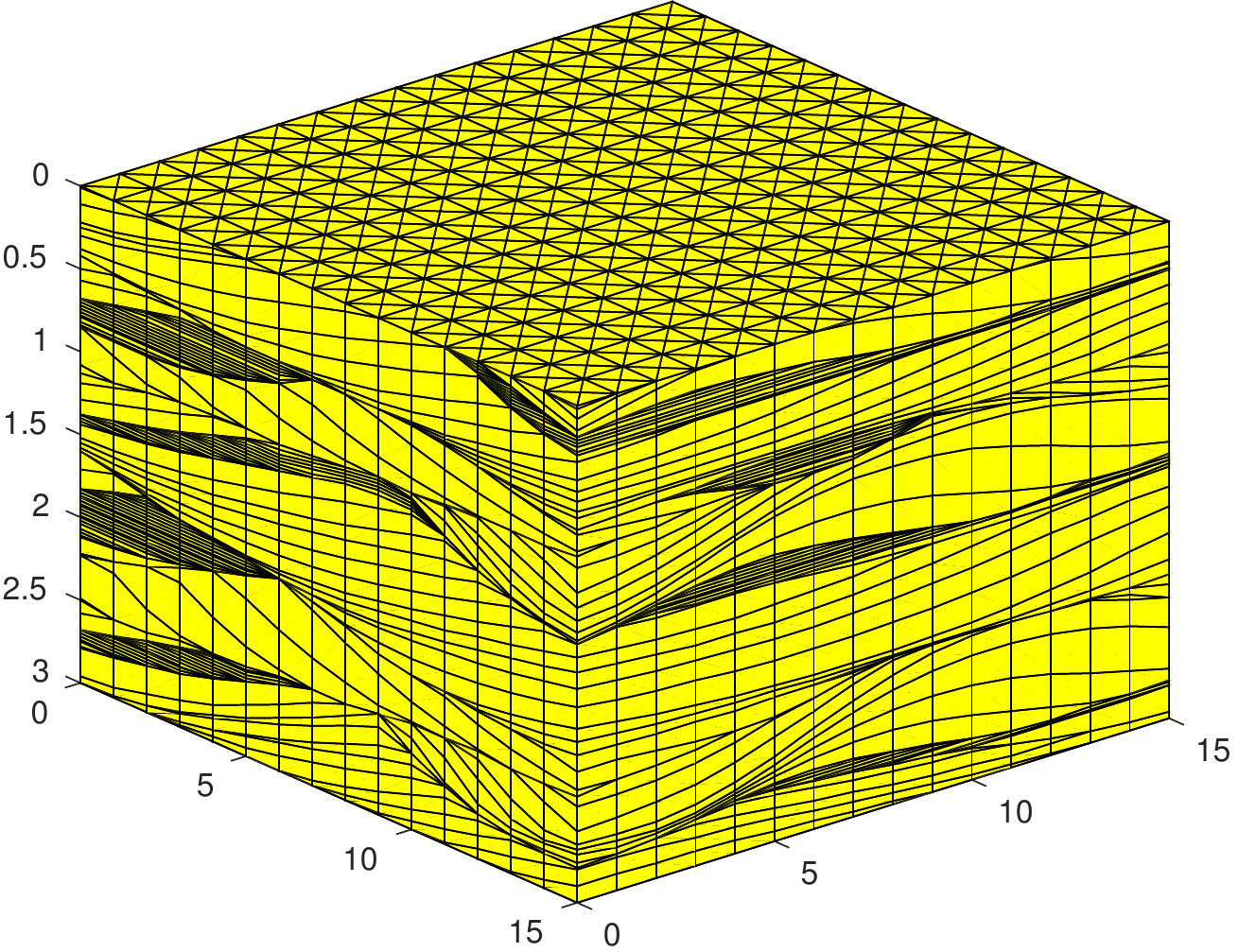}}; \\
      \node{\includegraphics[width=0.3\textwidth]{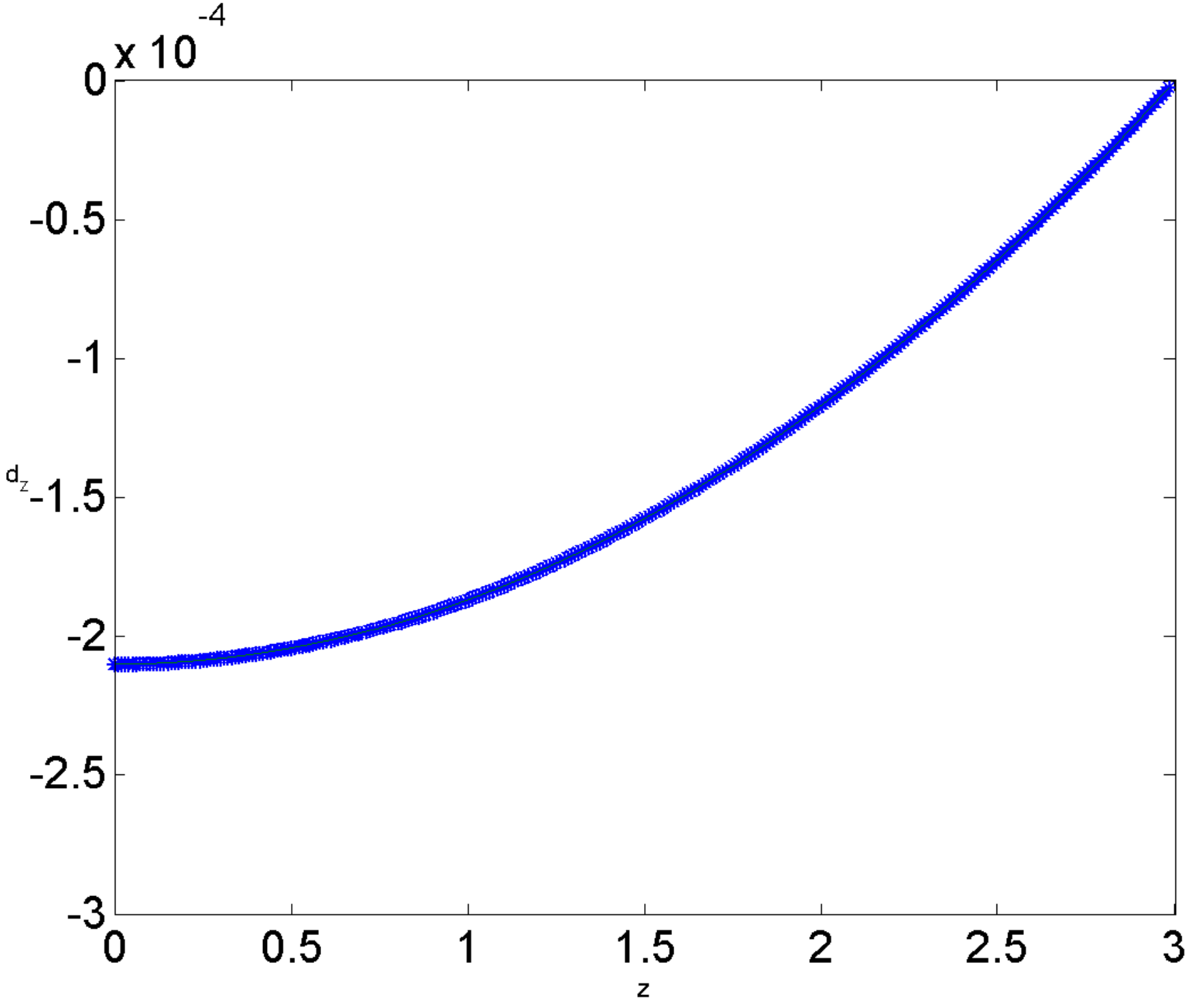}};&
      \node{\includegraphics[width=0.3\textwidth]{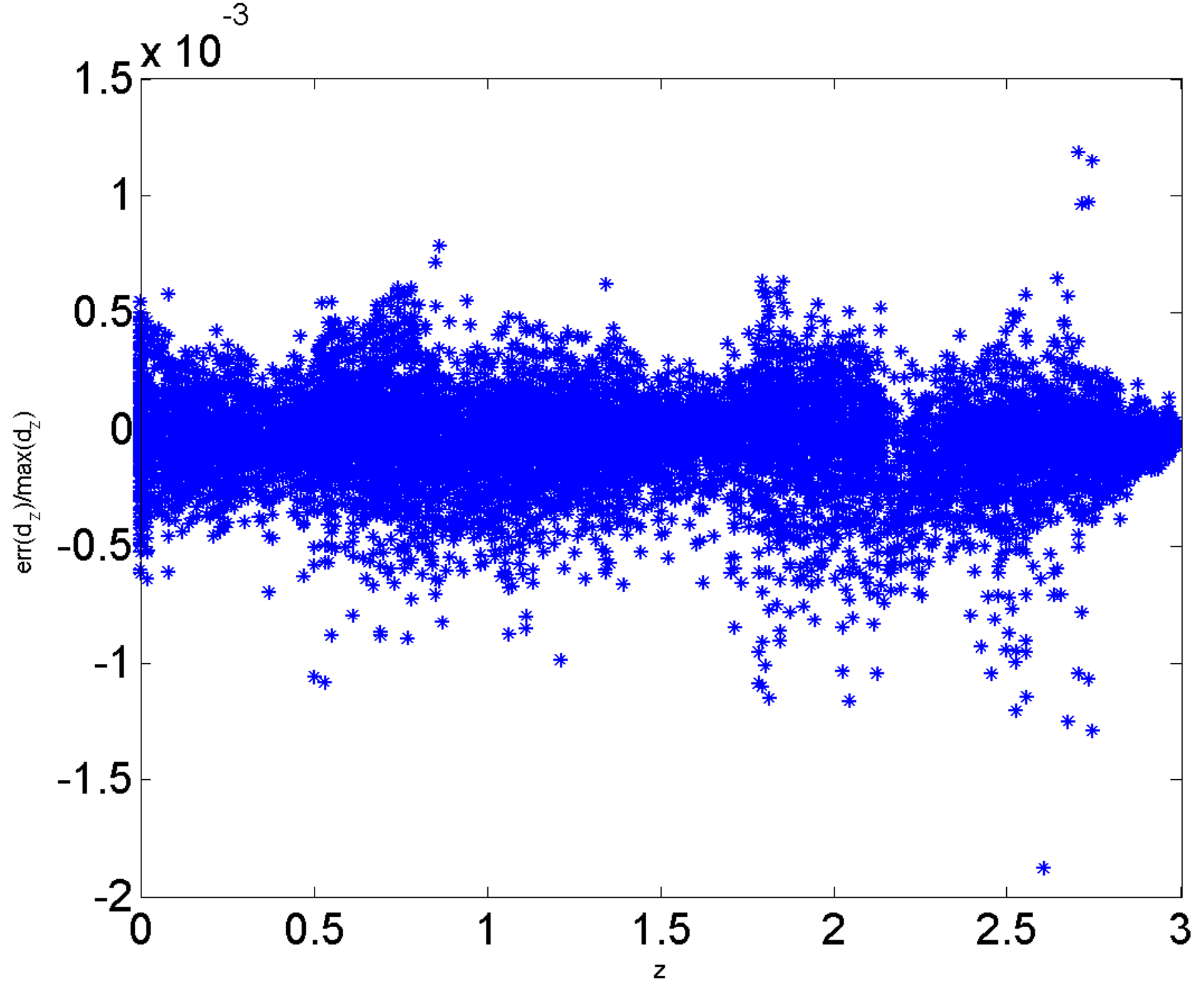}};&
      \node{\includegraphics[width=0.3\textwidth]{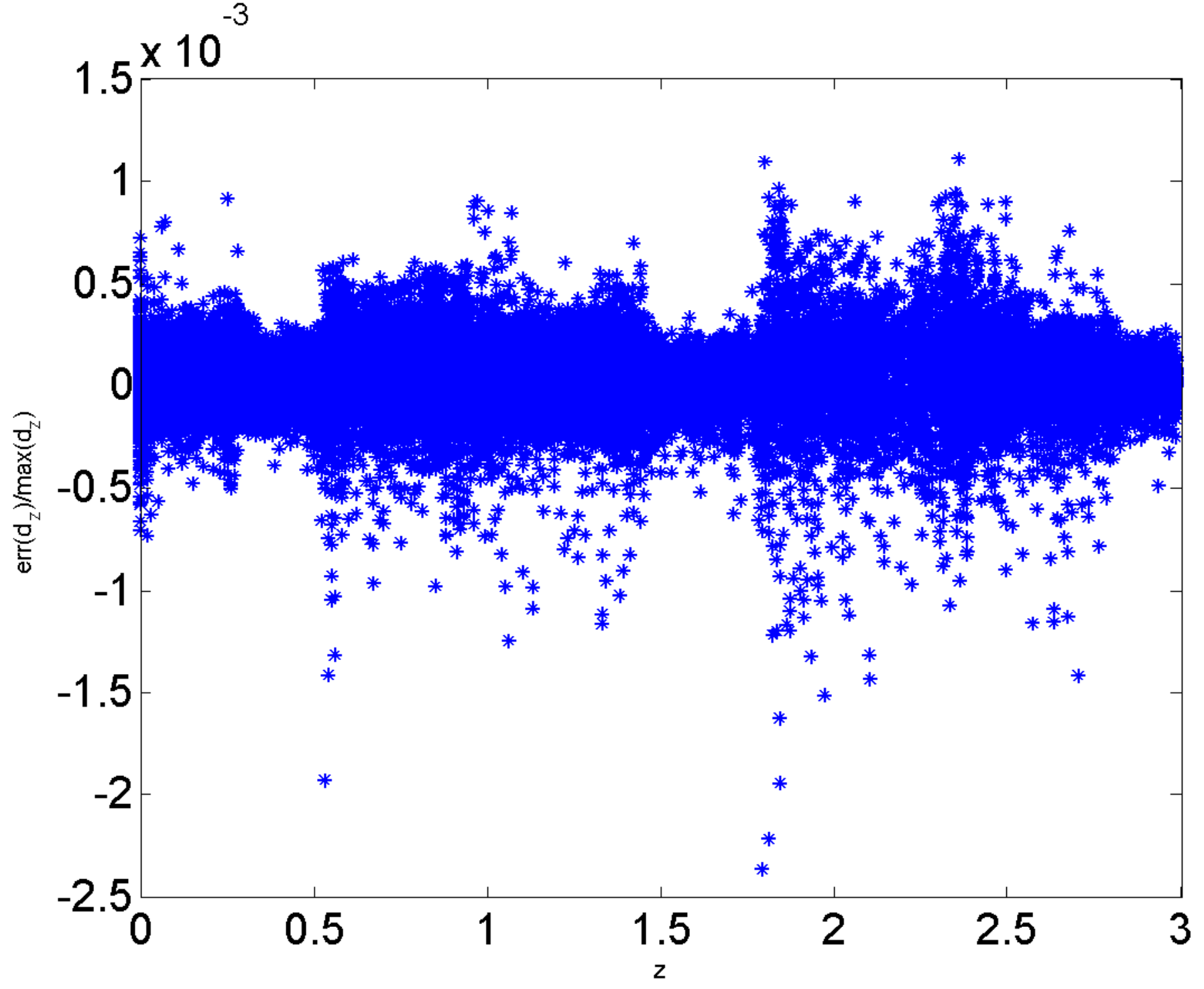}};\\
      \node{\includegraphics[width=0.3\textwidth]{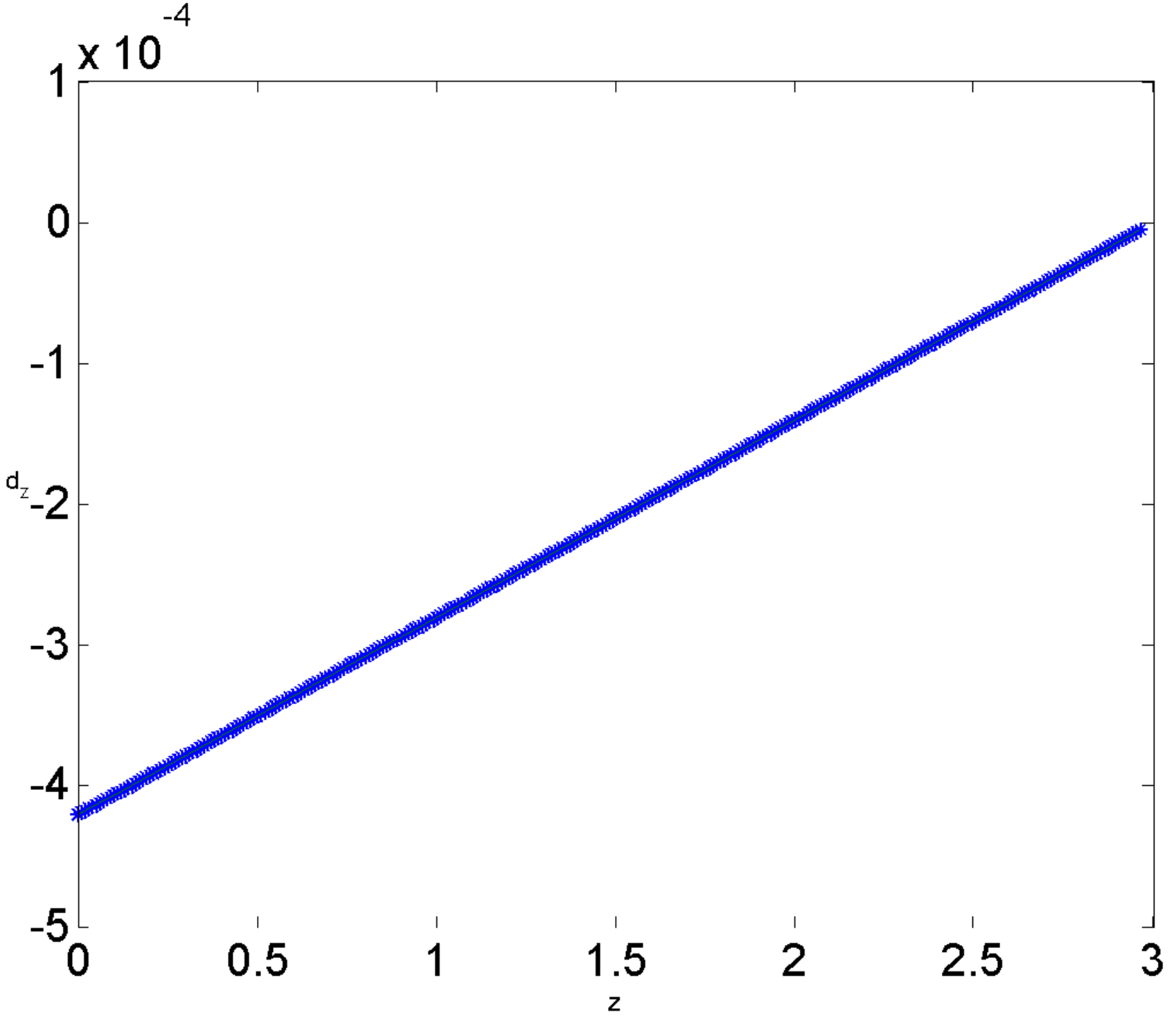}};&
      \node{\includegraphics[width=0.3\textwidth]{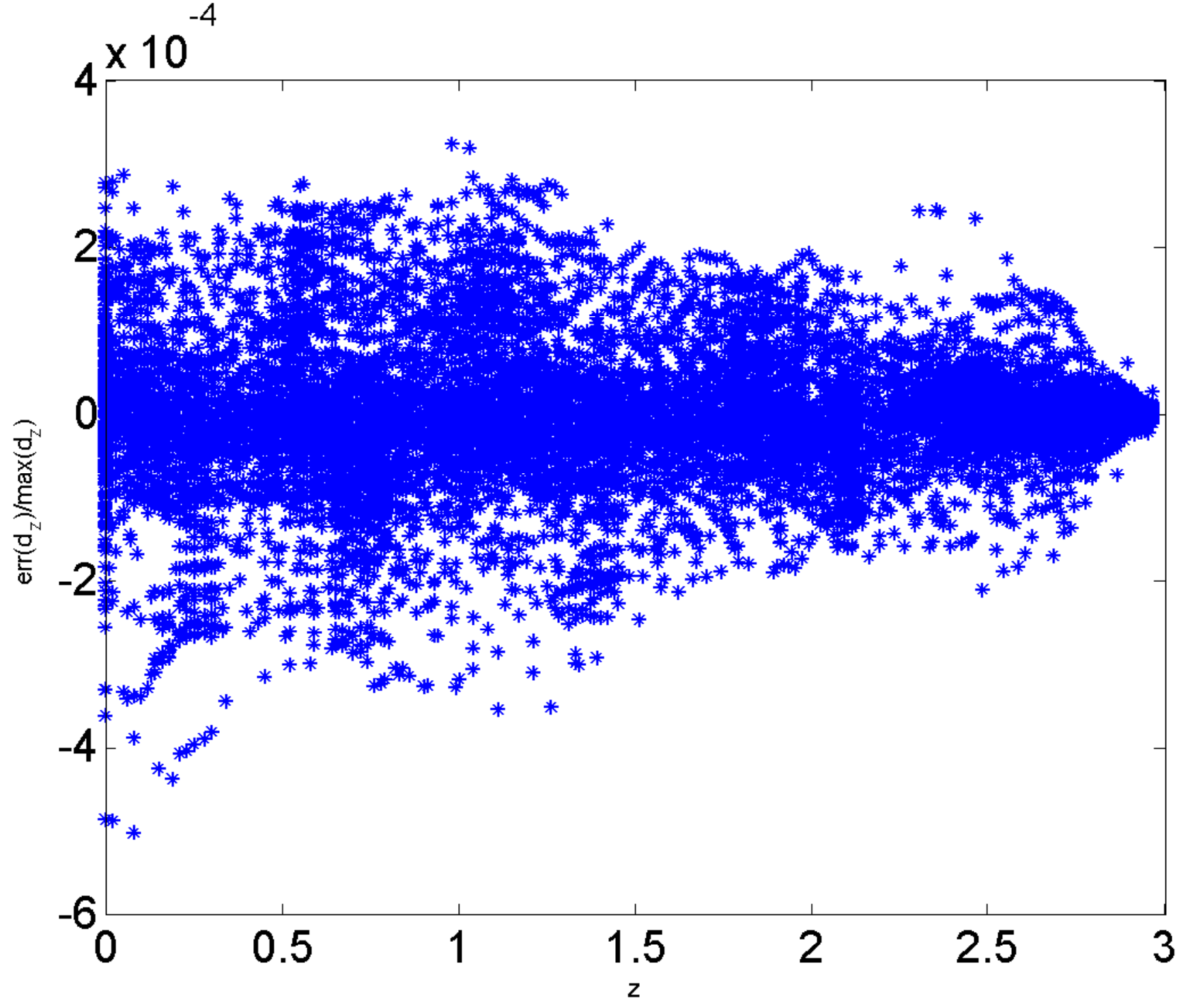}};&
      \node{\includegraphics[width=0.3\textwidth]{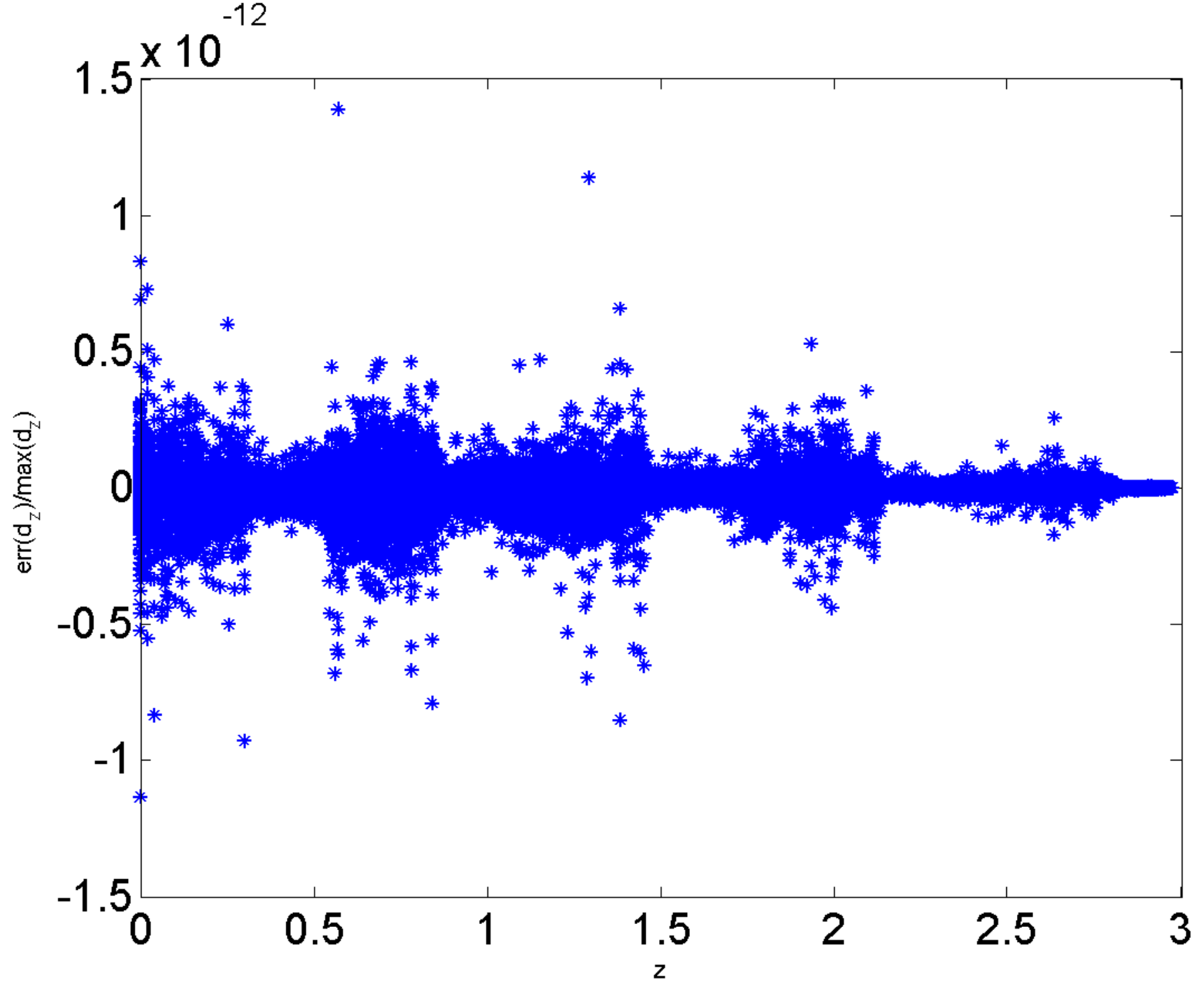}};\\
    };
    \node[anchor = south] at (notri.north) {Original corner-point grid};
    \node[anchor = south] at (tri.north) {Using triangulated faces};
  \end{tikzpicture}
  \caption{Effect of compression on the \sbed model. The first row shows a plot of
    the vertical deformation on the grid (left), the original grid (middle) and the
    same grid where the surface are triangulated (right). Two types of loads are
    considered: pure gravitational compression (second row), load at the top surface
    (third row). The first column shows the displacement obtained for each loading
    case, which is very close to the analytical solution. The remaining plots show
    the errors for the original cornerpoint grid with curved faces (middle column)
    and the triangulated grid with only planar faces (right column).}
  \label{fig:sbed}
\end{figure}
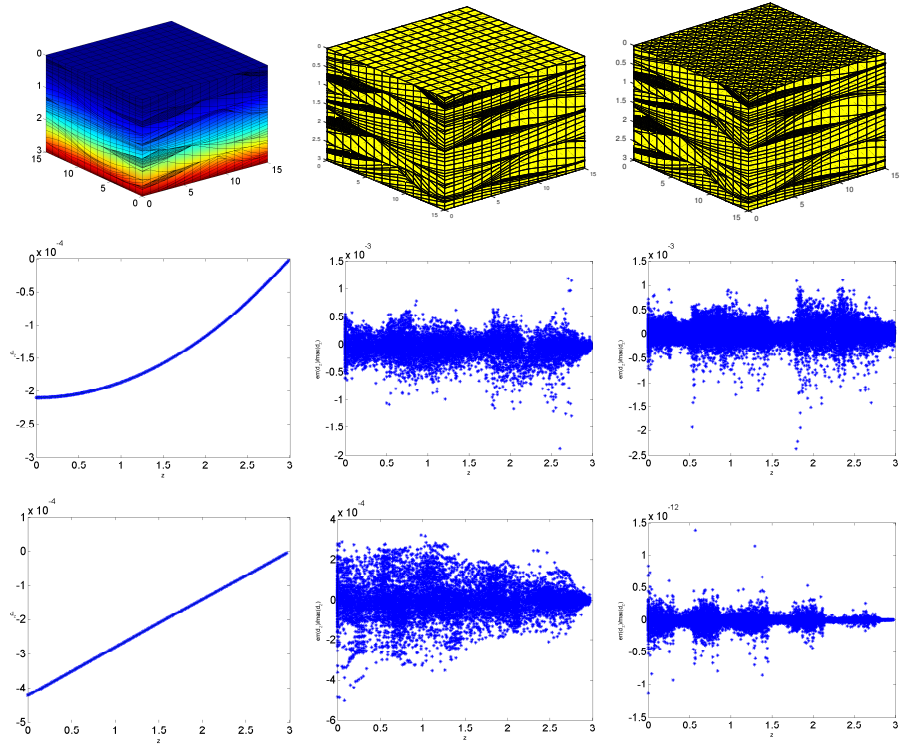

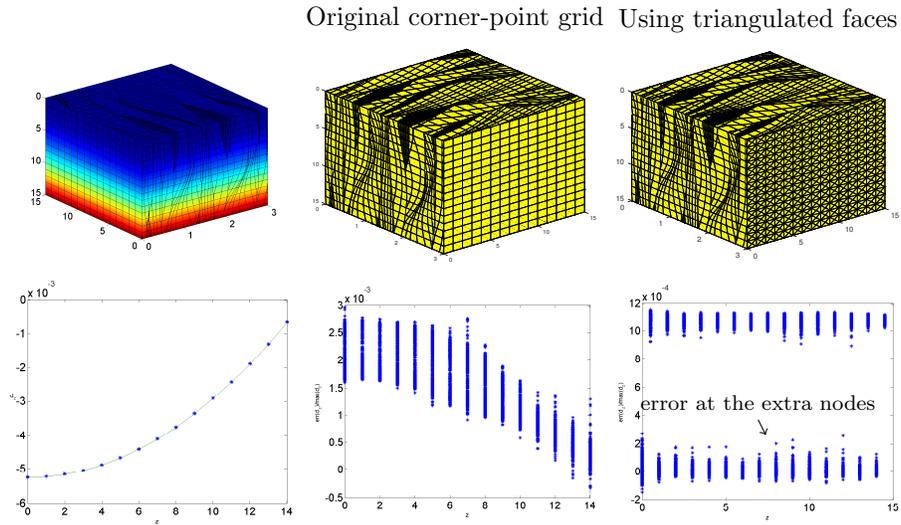
\begin{figure}
  \begin{tikzpicture}
    \matrix at (0,0) {
      \node{\includegraphics[width=0.3\textwidth]{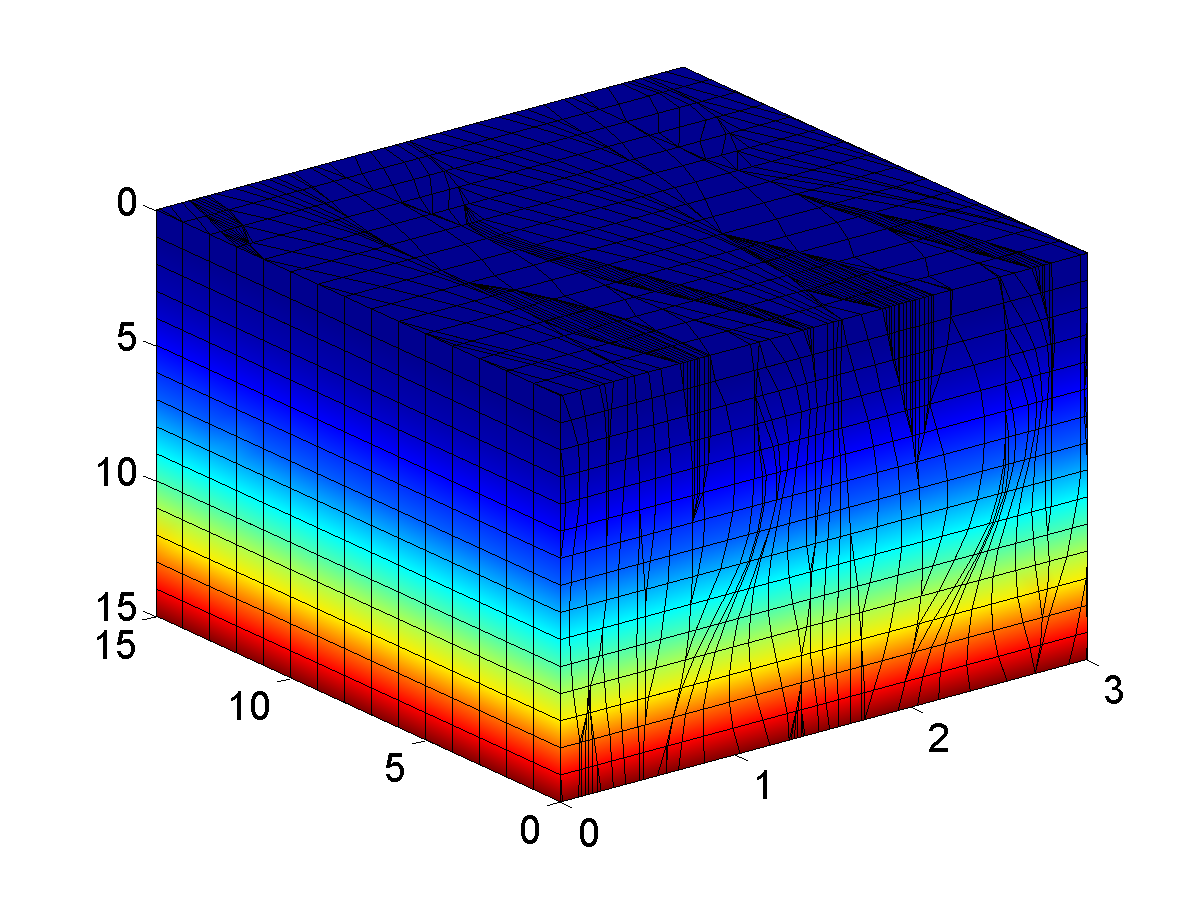}};
      &
      \node (notri) {\includegraphics[width=0.3\textwidth]{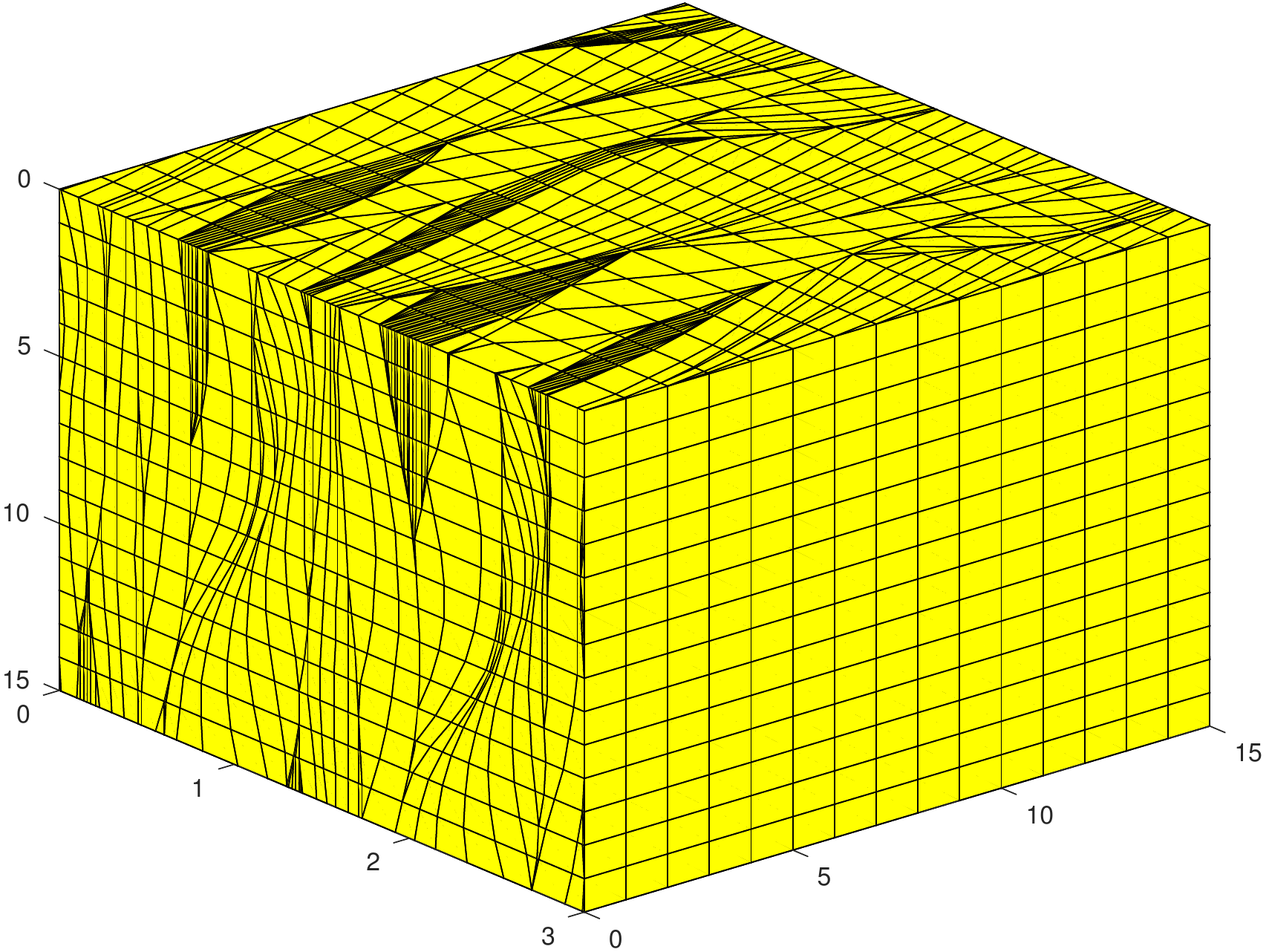}};
      &
      \node (tri)  {\includegraphics[width=0.3\textwidth]{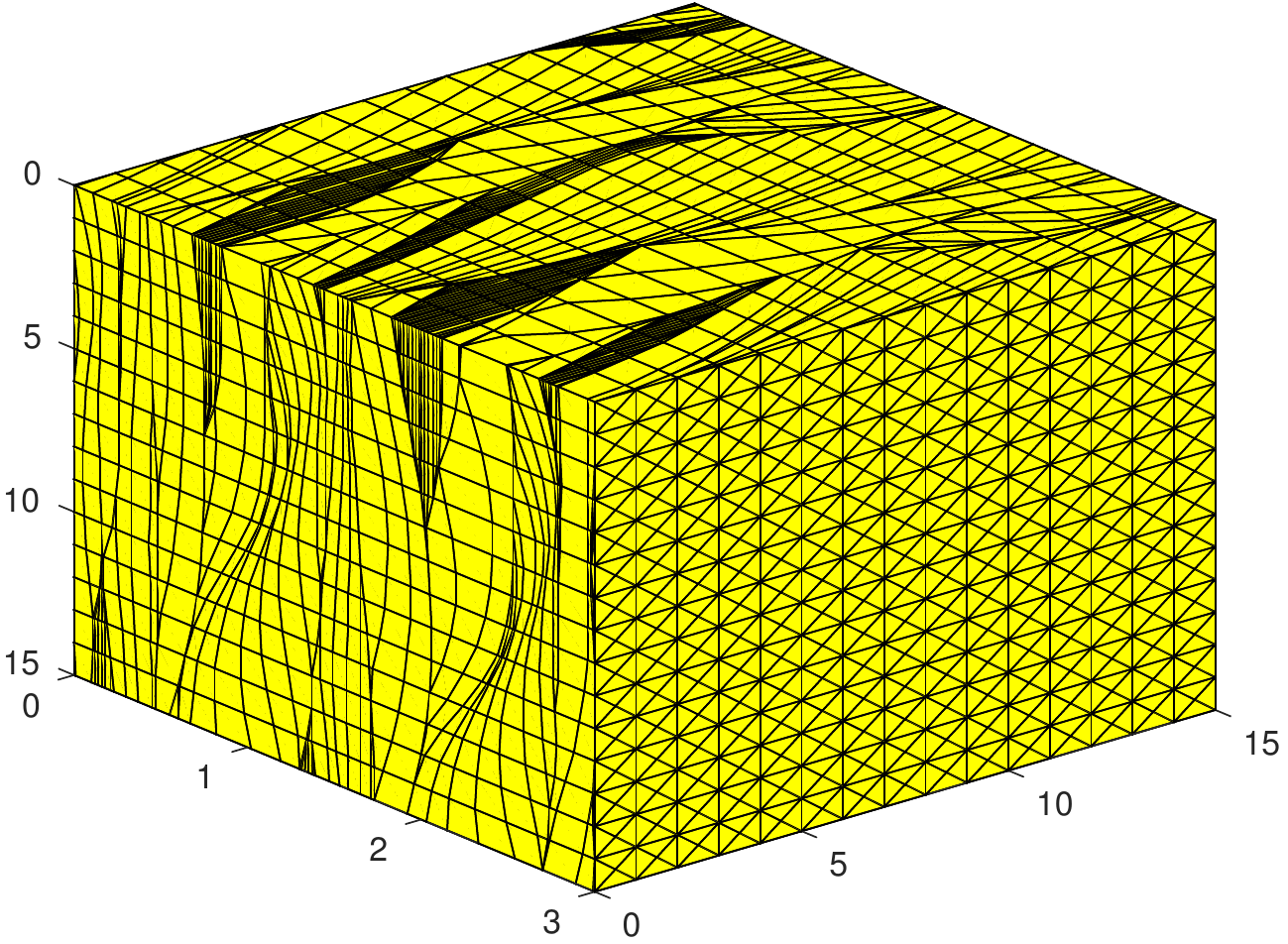}}; \\
      \node{\includegraphics[width=0.3\textwidth]{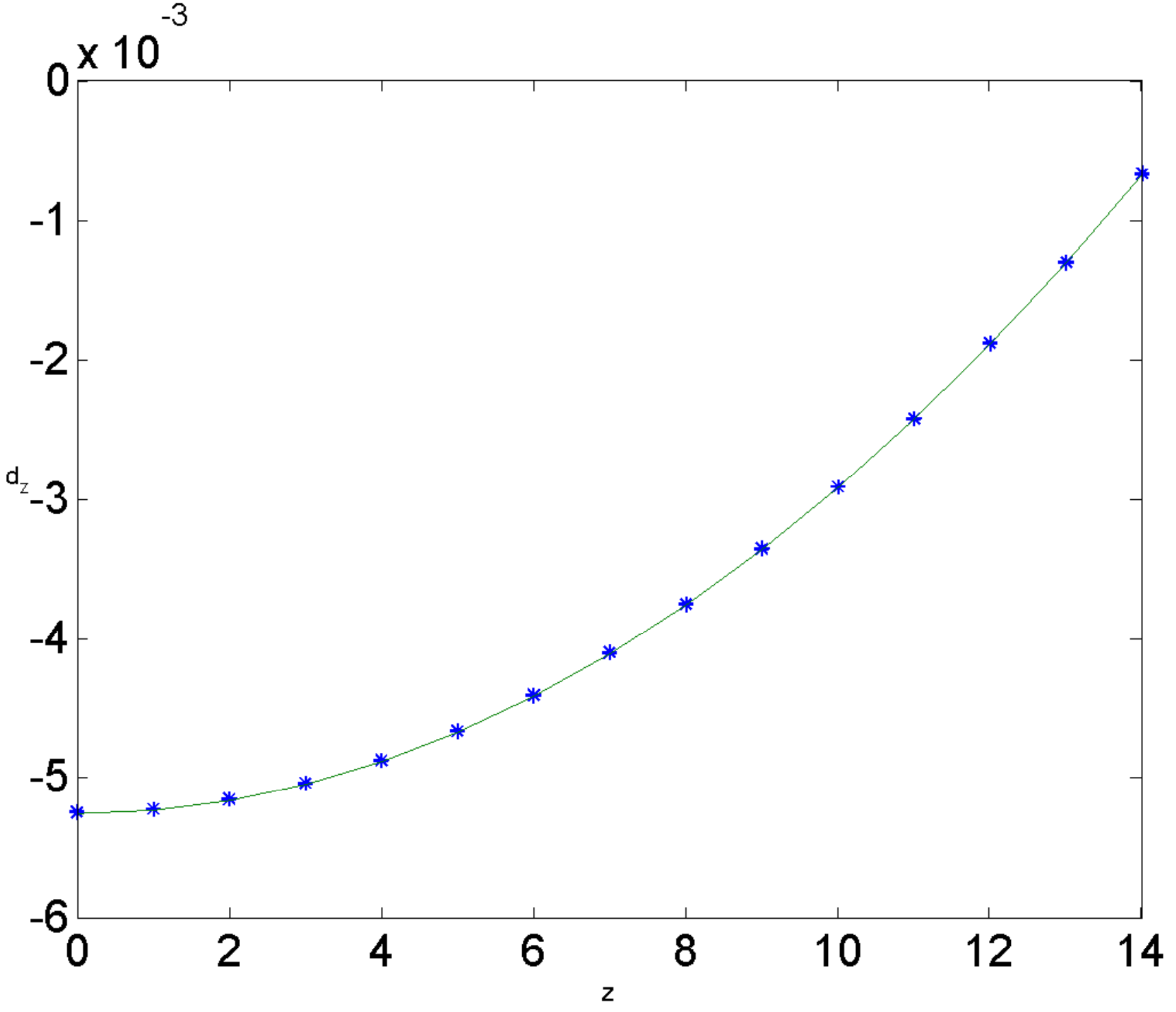}};&
      \node 
      {\includegraphics[width=0.3\textwidth]{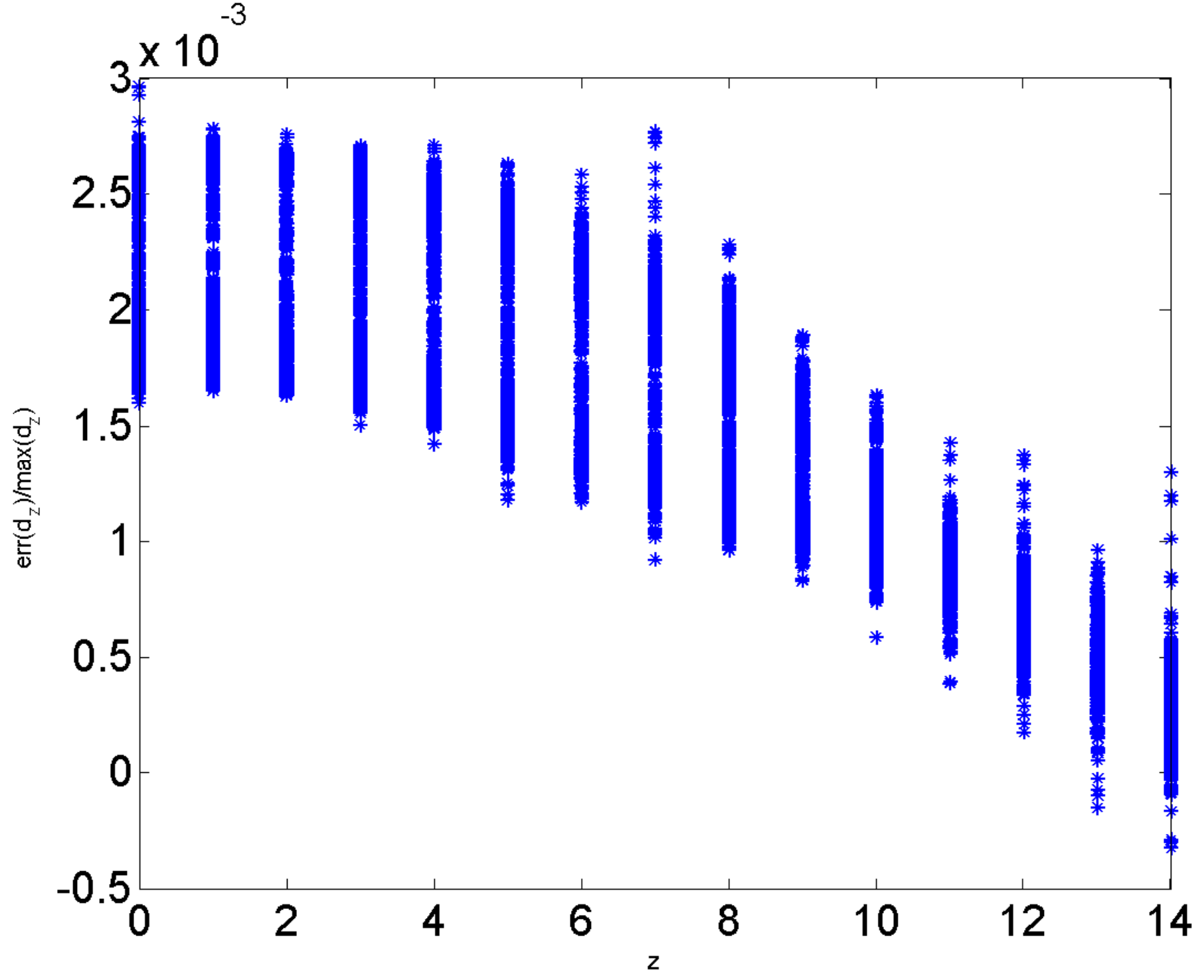}};&
      \node (splerr) 
      {\includegraphics[width=0.3\textwidth]{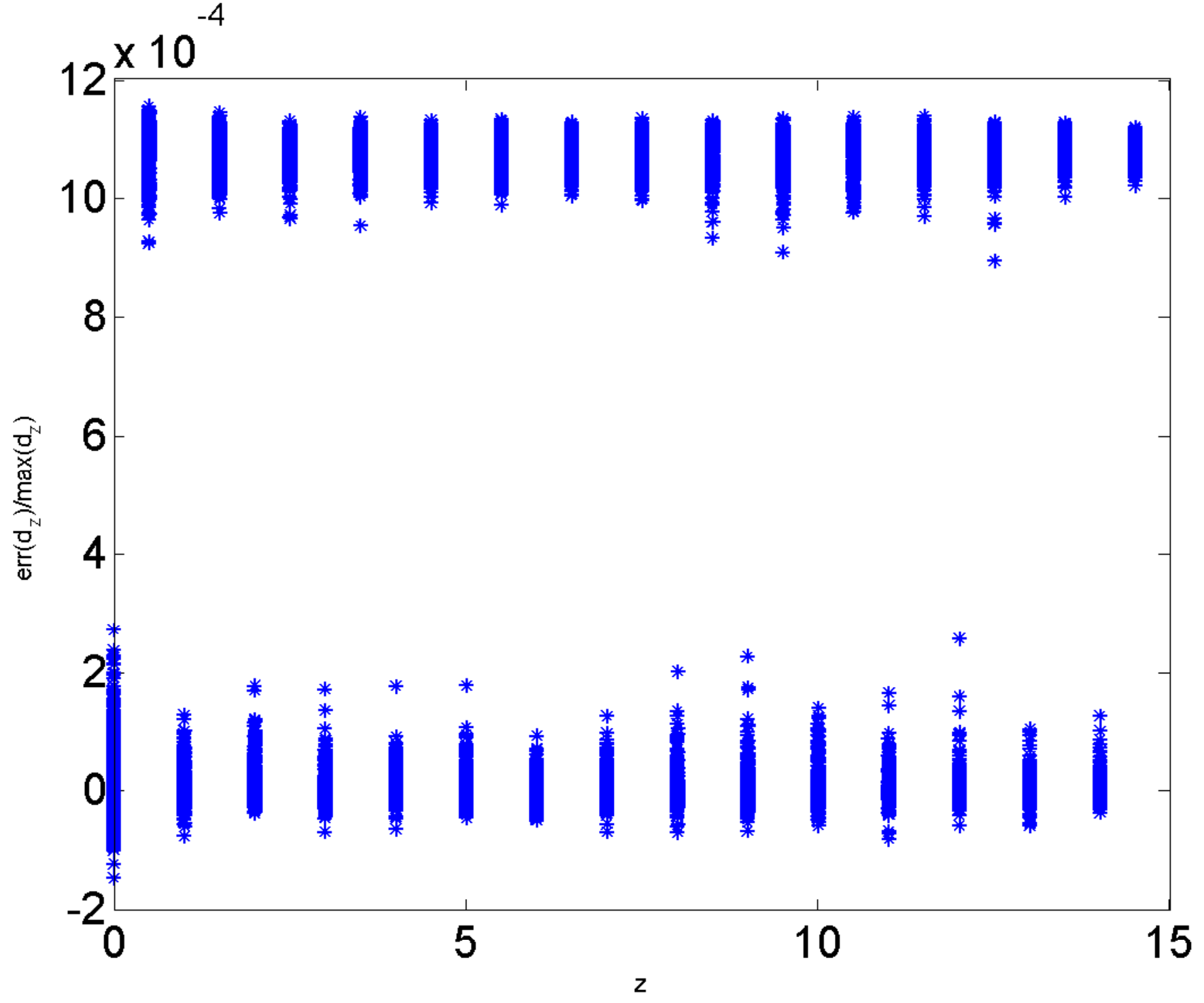}};\\};
    \node[anchor = south] at (notri.north) {Original corner-point grid};
    \node[anchor = south] at (tri.north) {Using triangulated faces};
    \node (exnodeleg) at (splerr.center) {\footnotesize error at the extra nodes};
    \draw [->] (exnodeleg.south) -- +(0.1, -0.2);
  \end{tikzpicture}
  \caption{Effect of compression on a fliped \sbed model. The first row shows a plot
    of the vertical deformation on the grid (left), the original grid (middle) and
    the same grid where the surface are triangulated (right). We consider only the
    case with gravitation load. The first column shows the displacement. The
    remaining plots show the errors for the original corner-point grid (middle
    column) and the triangulated grid (right column). On the plot at the lower right,
    we observe that the error splits clearly between the type of nodes, the extra
    face node at the bottom and the other at the top. }
  \label{fig:sbed_flip}
\end{figure}

In Figure \ref{fig:sbed}, we show the effect of compression with two types of load
given by a constant gravitational force and a constant load applied on the top
surface. For both loads, the analytical solutions can be computed and they are
respectively, quadratic and linear in $z$.  We consider both the original
corner-point grid and the triangulated grid. By triangulated grid, we mean a grid
where the faces ares triangulated, as we just explained. For all these cases, the VEM
method gives accurate results, given that we use the discrete gradient approach to
compute the load term. The other alternatives simply fail in this case, by errors
that are larger than the span of the exact solution. In the \sbed model, the pillars
are all vertical lines, which implies that the vertical faces are planar. For the
linear case corresponding to a constant load on the top surface, we see that the
triangulated version gives exact result, as predicted by the VEM theory, since in
this case all the surfaces are planar and the solution is linear. For the original
grid, we get an error due to the curved top and bottom faces in each cell. For the
pure gravitational case, both grids give comparable results. Thus, we can conclude
that in practice, it may not be worth triangulating the faces because it introduces
more degrees of freedom without significantly improving the accuracy of the
solution. We consider the case of a flipped model for Norne in Figure
\ref{fig:sbed_flip}. In this way, we can investigate the effect of having non planar
surface in the vertical direction. Typically, for the cells of the original
reservoir, we have $\frac{\Delta x}{\Delta z}\approx\frac{\Delta y}{\Delta z}\gg 1$
so that, by flipping the model, we can observe the consequence of inverting the
correlation between the aspect ratio and the direction of gravitation. The results
are similar. However the triangulated case which is exact for linear compression
highlights that the error of different types of nodes have different errors, see
explanation in the caption of Figure \ref{fig:sbed_flip}.

Besides features like layering and erosion, the Norne case introduces also fault
structures. Such grids are far from ideal for numerical calculation, but the VEM
method shows very robust behavior. In Figure \ref{fig:norne}, we look at the
difference between the original model and a model where all the pillars are
straightened up and made vertical. In this way, the curved sides in the 
vertical direction are eliminated. The analytical solution is unchanged
as we recall that the whole Norne model is anyway embedded in a regular prism.  The
results on Norne confirm those obtained for the \sbed model and show that effects of
curvature on the faces can be neglected. This indicates that, for many practical
applications, the VEM method can be used directly on the original grid of reservoirs
without deteriorating the accuracy of the results.

\clearpage

\section{Conclusion}

We have demonstrated how geomechanical calculations can be done directly on complex
geological models frequently encountered in reservoir modeling, by using the
flexibility of the VEM method which can handle general geometries. In this method,
the energy is not computed exactly for each basis element functions. We demonstrate
that this approximation can come at the cost of large errors for deformed grids, if
not care is taken when defining the approximate bilinear form. In particular we study
the effect of the load term calculation and show that, with stabilization terms and
load term calculations presented earlier in the literature, even simple 2D cases
fails severely when the aspect ratio is increased. We found that both the choices of
discretization and of the load term calculation are in combination responsible for
the failure. Using the exact equivalence with FEM on quadrilateral grid, we presented
a modification of the discretization that makes the method more robust in the 2D
case. In addition, we demonstrated that a calculation of the load in term of a
gradient of a potential was robust in 2D and the only approach which gave sufficient
accuracy in 3D. This holds in particular for grid cells that are outside the reach of
FEM, such as those containing hanging nodes. The VEM theory does not cover curved
faces, which are common in subsurface models. We saw that for our tests the error
associated with this feature was negligible comparable with other errors, with the
natural exception of the case when VEM gives the exact solution (linear
displacement).

\section*{acknowledgements}
This work has been partially funded by the Research Council of Norway through
  grants no. 215641 from the CLIMIT programme.


\end{document}